\documentclass{amsart}
\usepackage{amsfonts, amsbsy, amsmath, amssymb}

\ifx\fiverm\undefined 
  \newfont\fiverm{cmr5} 
\fi

\newtheorem{thm}{Theorem}[section]

\newtheorem{cla}[thm]{Claim}
\numberwithin{equation}{section}

\theoremstyle{definition}

\allowdisplaybreaks

\newcommand{\f}{\Bbb F}

\begin{document}

\title{On the Tu-Zeng Permutation Trinomial of Type $(1/4,3/4)$}

\author[Xiang-dong Hou]{Xiang-dong Hou}
\address{Department of Mathematics and Statistics,
University of South Florida, Tampa, FL 33620}
\email{xhou@usf.edu}



\begin{abstract}
Let $q$ be a power of $2$. Recently, Tu and Zeng considered trinomials of the form $f(X)=X+aX^{(1/4)q^2(q-1)}+bX^{(3/4)q^2(q-1)}$, where $a,b\in\Bbb F_{q^2}^*$. They proved that $f$ is a permutation polynomial of $\Bbb F_{q^2}$ if $b=a^{2-q}$ and $X^3+X+a^{-1-q}$ has no root in $\Bbb F_q$. In this paper, we show that the above sufficient condition is also necessary.
\end{abstract}

\maketitle

\section{Introduction}

Let $\f_q$ denote the finite field with $q$ elements. A polynomial $f\in\f_q[X]$ is called a permutation polynomial (PP) of $\f_q$ if it induces a permutation of $\f_q$. Consider a particular type of polynomials $f\in\f_q[X]$, e.g., binomials, trinomials, etc., to determine the PPs of the given type is to find criteria (in terms of the coefficients of $f$) for the equation $f(x)=y$ to have at least (or at most) one solution $x\in\f_q$ for each $y\in\f_q$. Such questions are usually difficult. Much of the research in this direction have been focusing on trinomials with Niho exponents. These are polynomials of the form
\begin{equation}\label{1.1}
f_{(r,s_1,s_2,),a,b}(X)=X^r(1+aX^{s_1(q-1)}+bX^{s_2(q-1)})\in\f_{q^2}[X],
\end{equation}
where $r,s_1,s_2$ are positive integers and $a,b\in\f_{q^2}$, and the objective is to determine the conditions on $a$ and $b$ that are necessary and sufficient for $f_{(r,s_1,s_2),a,b}$ to be a PP of $\f_{q^2}$. The integers $s_1$ and $s_2$ can be treated as elements of $\Bbb Z/(q+1)\Bbb Z$; for example, for even $q$, $(s_1,s_2)=(1/4,3/4)$ is meaningful. The question has been solved for the following cases of parameters $(r,s_1,s_2)$:
\begin{itemize}
\item $(r,s_1,s_2)=(1,1,2)$ \cite{Hou-FFA-2015b}.

\item $(r,s_1,s_2)=(1,-1/2,1/2)$, $\text{char}\,\f_q=2$ \cite{Tu-Zeng-FFA-2018}.

\item $(r,s_1,s_2)=(1,-1,2)$, $\text{char}\,\f_q=2,3$ \cite{Bartoli-FFA-2018, Hou-arXiv1803.04071, Hou-Tu-Zeng-arXiv1811.11949, Tu-Zeng-Li-Helleseth-FFA-2015}.
\end{itemize}
The purpose of the present paper is to add the case $(r,s_1,s_2)=(1,1/4,3/4)$, $\text{char}\,\f_q=2$, to the above short list. Let $q$ be even. In a recent paper \cite{Tu-Zeng-FFA-2018}, Tu and Zeng proved that $f_{(1,1/4,3/4),a,b}(X)=X(1+aX^{(q-1)/4}+bX^{3(q-1)/4})$ $(a,b\in\f_{q^2}^*$) is a PP of $\f_{q^2}$ if $b=a^{2-q}$ and $X^3+X+a^{-1-q}$ has no root in $\f_q$. Note that $f_{(1,1/4,3/4),a,b}(X^4)=X^4(1+aX^{q-1}+bX^{3(q-1)})=f_{(4,1,3),a,b}(X)$, and by a suitable substitution $X\mapsto uX$, $u\in\f_{q^2}^*$, we may assume that $a\in\f_q^*$. Therefore, the result of Tu and Zeng can be stated as follows: Let $f(X)=X^4(1+aX^{q-1}+bX^{3(q-1)})$, where $a\in\f_q^*$ and $b\in\f_{q^2}^*$. Then $f$ is a PP of $\f_{q^2}$ if $a=b$ and $X^3+X+a^{-1}$ has no root in $\f_q$. We show that the sufficient condition here is also necessary:

\begin{thm}\label{T1.1}
Let $q$ be even and $f(X)=X^4(1+aX^{q-1}+bX^{3(q-1)})$, where $a\in\f_q^*$ and $b\in\f_{q^2}^*$. Then $f$ is a PP of $\f_{q^2}$ if and only if $a=b$ and $X^3+X+a^{-1}$ has no root in $\f_q$.
\end{thm}

The proof of Theorem~\ref{T1.1} follows a strategy similar to that of \cite{Hou-arXiv1803.04071}. By a well-known folklore, we see that $f$ is a PP of $\f_{q^2}$ if and only if, essentially, a certain quartic equation in $X$ with coefficients in $\f_q(Y)$ has a unique solution $X=x\in\f_q$ for every $Y=y\in\f_q$. This fact, combined with a theorem by Leonard and Williams on the factorization of quartic polynomials over $\f_q$ and a theorem by Hou and Iezzi on composition of rational functions over $\f_q$, implies a compositional factorization of a certain rational function in $\f_q(Y)$. Comparison of the coefficients in the factorization gives several polynomial equations involving $a$ and $b$. A careful analysis of these equations, mostly done through computerized symbolic computations, concludes that the solutions $(a,b)$ of these equations, filtered by some additional requirements, are  precisely those described in Theorem~\ref{T1.1}.

Theorem~\ref{T1.1} is proved in three steps, in the remaining three sections, respectively, each with a specific goal. In Section 2, we show that if $f$ is a PP of $\f_{q^2}$, then a certain rational function in $X$ and $Y$ has a unique zero $(x,y)\in\f_q^2$ for every $y\in\f_q$. In Section 3, we prove that if $f$ is a PP of $\f_{q^2}$, then $b\in\f_q^*$. This is the most difficult part of the proof. In Section 4, under the assumption that $b\in\f_q^*$, we conclude that $f$ is a PP of $\f_{q^2}$ if and only if $a=b$ and $X^3+X+a^{-1}$ has no root in $\f_q$.  

The following two results will be used in the proof of Theorem~\ref{T1.1}:

\begin{thm}[Leonard and Williams~\cite{Leonard-Williams-PAMS-1972}, Williams \cite{Williams-JNT-1975}] \label{LW} 
 Let $q$ be even and let and $\alpha_0,\alpha_1,\alpha_2\in\f_q$ be such that $\alpha_0\alpha_1\ne 0$.
\begin{itemize}
\item[(i)]  $X^4+\alpha_2 X^2+\alpha_1X+\alpha_0$ has a unique root in $\f_q$ if and only if $X^3+\alpha_2X+\alpha_1$ is irreducible over $\f_q$.

\item[(ii)] $X^3+\alpha_2X+\alpha_1$ is irreducible over $\f_q$ if and only if $\text{\rm Tr}_{q/2}(1+\alpha_2^3/\alpha_1^2)=0$ and $X^6+\alpha_1X^3+\alpha_2^3$ has no root in $\f_{q^2}$.
\end{itemize}
\end{thm}

\begin{thm}[Hou and Iezzi \cite{Hou-Iezzi-ppt}]\label{HI}
Let $F(X),G(X)\in\f_q(X)\setminus\f_q$ be such that $\deg F=d$ and $\deg G=\delta$. If there is a constant $0<\epsilon\le 1$ such that
\begin{equation}\label{3.1}
|\{(x,y)\in\f_q\times\f_q:F(x)=G(y)\}|\ge q\Bigl(\Bigl\lfloor\frac\delta 2\Bigr\rfloor+\epsilon\Bigr),
\end{equation}
and $q\ge(d+\delta)^4/\epsilon^2$, then $F(X)=G(H(X))$ for some $H\in\f_q(X)$.
\end{thm}

In our notation, the resultant of two polynomials $P_1(X)$ and $P_2(X)$ is denoted by $\text{Res}(P_1,P_2;X)$. If $P\in\f_q[X]$ and $u\in\f_q$ is an unknown element, $P(u)$ is frequently treated as a polynomial in $u$ rather than an element of $\f_q$. Therefore the meaning of $\text{Res}(P_1(u),P_2(u);u)$ is $\text{Res}(P_1(X),P_2(X);X)|_{X=u}$.

\section{Initial Approach}

From now on, $q$ is even and
\begin{equation}\label{2.1}
f(X)=X^4(1+aX^{q-1}+bX^{3(q-1)}),
\end{equation}
where $a\in\f_q^*$ and $b\in\f_{q^2}^*$.  We assume that $q\ge 2^{13}$. (For $q\le 2^{12}$, Theorem~\ref{T1.1} can be verified with a computer.) It is well known \cite{Park-Lee-BAMS-2001, Wang-LNCS-2007, Zieve-PAMS-2009} that $f$ is a PP of $\f_{q^2}^*$ if and only if $f_1(X):=X^4(1+aX+bX^3)^{q-1}$ permutes $\mu_{q+1}:=\{x\in\f_{q^2}:x^{q+1}=1\}$. For $x\in\mu_{q+1}$ with $1+ax+bx^3\ne 0$, we have $f_1(x)=g(x)$, where
\begin{equation}\label{2.2}
g(X)=\frac{X(b^q+aX^2+X^3)}{1+aX+bX^3}.
\end{equation}
Therefore, $f$ is a PP of $\f_{q^2}$ if and only if $1+aX+bX^3$ has no root in $\mu_{q+1}$ and $g(X)$ permutes $\mu_{q+1}$.

Assume that $1+aX+bX^3$ has no root in $\mu_{q+1}$, in particular, $1+a+b\ne0$. Choose $z\in\f_{q^2}$ such that $\text{Tr}_{q^2/q}(z)=1$ and let $k=\text{N}_{q^2/q}(z)$. Then $z^2+z+k=0$ and $\text{Tr}_{q/2}(k)=1$. Let $\phi(X)=(X+z^q)/(X+z)$, which maps $\f_q\cup\{\infty\}$ to $\mu_{q+1}$ bijectively, and let $\psi(X)=(1+a+b)^{q-1}\phi(X)$. We have the following diagram:
\[
\renewcommand{\arraystretch}{1.5}
\begin{array}{ccccccc}
\f_q\cup\{\infty\}&\overset{\phi}\longrightarrow &\mu_{q+1} &\overset{g}\longrightarrow &\mu_{q+1} & \overset{\psi^{-1}}\longrightarrow &\f_q\cup\{\infty\}\\
\infty &\longmapsto & 1 &\longmapsto &(1+a+b)^{q-1}&\longmapsto &\infty
\end{array}
\]
where $\psi^{-1}$ denotes the compositional inverse of $\psi$. Therefore, $g(X)$ permutes $\mu_{q+1}$ if and only if $\psi^{-1}\circ g\circ \phi$ permutes $\f_q$, that is, if and only if for each $y\in\f_q$, there is a unique $x\in\f_q$ such that $(g\circ\phi)(x)=\psi(y)$, i.e.,
\begin{equation}\label{2.3}
g\Bigl(\frac {x+z+1}{x+z}\Bigr)=(1+a+b)^{q-1}\frac{y+z+1}{y+z}.
\end{equation}

\section{Necessity that $b\in\f_q^*$}

In this section we prove the following claim:

\begin{cla}\label{C3.1}
If $f$ is a PP of $\f_{q^2}$, then $b\in\f_q^*$.
\end{cla}

The proof of this moderate claim turns out to be quite involved. For the rest of this section, assume that $f$ is a PP of $\f_{q^2}$ and assume to the contrary that $b\notin \f_q$.

\subsection{From a quartic equation to a system of equations} \label{s3.1}

Let $b_1=b+b^q$ and $z=b/b_1$. Then $b=b_1z$, $\text{Tr}_{q^2/q}(z)=1$, and $k:=\text{N}_{q^2/q}(z)=b^{1+q}/b_1^2$. In \eqref{2.3}, write 
\begin{equation}\label{2.4.0}
g\Bigl(\frac{X+z+1}{X+z}\Bigr)=\frac{A(X)}{B(X)},
\end{equation}
where
\begin{align*}
A(X)=\,& 1+az+az^2+az^3+b_1z^3+z^4+az^4+b_1z^5 +(a+az^2+b_1z^2+b_1z^3 )X\cr
& +(a+az+b_1z+b_1z^2 )X^2 +(a+b_1+b_1z )+X^3 +(1+a+b_1+b_1z  )X^4, \cr
B(X)=\,& b_1z^2+az^3+b_1z^3 +z^4+az^4+b_1z^4+b_1z^5 + (b_1z+az^2+b_1z^3 )X\cr
& +(az+b_1z+b_1z^2 )X^2 +(a+b_1z )X^3 +(1+a+b_1z )X^4.
\end{align*}
Since $1+aX+bX^3$ has no root in $\mu_{q+1}$, $B(X)$ has no root in $\f_q$. Now \eqref{2.3} becomes
\begin{equation}\label{2.4}
\frac{A(x)}{B(x)}=\frac{1+a+b_1(z+1)}{1+a+b_1z}\cdot\frac{y+z+1}{y+z}.
\end{equation}
Using the relation $z^2+z+k=0$, we write \eqref{2.4} as
\begin{equation}\label{2.5}
C_4x^4+C_3(y)x^3+C_2(y)x^2+C_1(y)x+C_0(y)=0,
\end{equation}
where
\begin{align*}
C_4=\,& a^2+a b_1+b_1^2 k+b_1+1\ (\ne 0), \cr
C_3(Y)=\,& a^2+a b_1+a+b_1^2 k+b_1 Y, \cr
C_2(Y)=\,& b_1^2 k+b_1 k + (a^2+a+b_1^2 k)Y, \cr
C_1(Y)=\,& a^2 k+a
b_1 k+a k+b_1^2 k^2+b_1^2 k+b_1 k + (a^2+a+b_1^2 k+b_1 k)Y, \cr
C_0(Y)=\,& a^2 k^2+a
b_1 k^2+a k+b_1^2 k^3+b_1^2 k^2+k^2+k\cr
&+ (a^2 k+a k+a+b_1^2 k^2+b_1 k+1)Y.
\end{align*}
Recall from Section 2 that for each $y\in\f_q$, there is a unique $x\in\f_q$ satisfying \eqref{2.5}. In the above $C_4\ne 0$ since otherwise, $\text{Tr}_{q/2}(k)=\text{Tr}_{q/2}(b_1^{-1}+b_1^{-2}+a/b_1+(a/b_1)^2)=0$, which is a contradiction. Write $C_i=C_i(Y)$ and $c_i=C_i(y)$, $0\le i\le 4$. Squaring \eqref{2.5} and replacing $x$ with $x^{1/2}$ give 
\begin{equation}\label{2.6}
c_4^2x^4+c_3^2x^3+c_2^2x^2+c_1^2x+c_0^2=0.
\end{equation}
Clearly, $C_3\ne 0$. We also have $C_1^2C_4+C_1C_2C_3+C_0C_3^2\ne 0$; this is equivalent to saying that $E_0,\dots,E_3$, given in \eqref{E0} -- \eqref{E3}, are not all $0$. In fact, at the end of this subsection, we show that $E_2$ and $E_3$ are not both $0$. 
Let $y\in\f_q$ be such that $c_3(c_1^2c_4+c_1c_2c_3+c_0c_3^2)\ne 0$. Replacing $x$ be $x+c_1/c_3$ in \eqref{2.6} gives
\[
c_3^4c_4^2x^4+c_3^6x^3+c_3^4(c_1c_3+c_2^2)x^2+(c_1^2c_4+c_1c_2c_3+c_0c_3^2)^2=0.
\]
Further replacing $x$ by $1/x$ in the above gives
\begin{equation}\label{2.7}
x^4+\alpha_2x^2+\alpha_1x+\alpha_0=0,
\end{equation}
where
\begin{align*}
\alpha_0=\,&\frac{c_3^4c_4^2}{(c_1^2c_4+c_1c_2c_3+c_0c_3^2)^2},\cr
\alpha_1=\,&\frac{c_3^6}{(c_1^2c_4+c_1c_2c_3+c_0c_3^2)^2},\cr
\alpha_2=\,&\frac{c_3^4(c_1c_3+c_2^2)}{(c_1^2c_4+c_1c_2c_3+c_0c_3^2)^2}.
\end{align*}
Therefore \eqref{2.7} has a unique solution $x$ in $\f_q$. By Theorem~\ref{LW}, we have
\[
\text{Tr}_{q/2}\Bigl(1+\frac{\alpha_2^3}{\alpha_1^2}\Bigr)=0,
\]
i.e.,
\[
\text{Tr}_{q/2}\Bigl(1+\frac{(c_1c_3+c_2^2)^3}{(c_1^2c_4+c_1c_2c_3+c_0c_3^2)^2}\Bigr)=0.
\]
To recap, for each $y\in\f_q$ such that $c_3(c_1^2c_4+c_1c_2c_3+c_0c_3^2)\ne 0$, there are two $x\in\f_q$ such that 
\[
1+\frac{(c_1c_3+c_2^2)^3}{(c_1^2c_4+c_1c_2c_3+c_0c_3^2)^2}=x+x^2.
\]
Let $F(Y)=1+(C_1C_3+C_2^2)^3/(C_1^2C_4+C_1C_2C_3+C_0C_3^2)^2$ and $G(Y)=Y+Y^2$, where $\deg F=d\le 6$ and $\deg G=\delta=2$. Choose $\epsilon=2^{-1/2}$. Then we have
\[
|\{(y,x):\in\f_q\times\f_q:F(y)=G(x)\}|\ge (q-4)\cdot 2=2q-8>q(1+\epsilon)
\]
and 
\[
q\ge 2^{13}=(6+2)^4/\epsilon^2.
\]
By Theorem~\ref{HI}, there exists $H\in\f_q(Y)$ such that $F(Y)=G(H(Y))$, i.e.,
\begin{equation}\label{2.8}
1+\frac{(C_1C_3+C_2^2)^3}{(C_1^2C_4+C_1C_2C_3+C_0C_3^2)^2}=H(Y)+H(Y)^2.
\end{equation}
It follows that 
\[
H(Y)=\frac{D(Y)}{C_1^2C_4+C_1C_2C_3+C_0C_3^2}
\]
for some $D(Y)\in\f_q[Y]$ and \eqref{2.8} becomes
\begin{equation}\label{2.9}
(C_1^2C_4+C_1C_2C_3+C_0C_3^2)^2+(C_1C_3+C_2^2)^3=D(D+C_1^2C_4+C_1C_2C_3+C_0C_3^2).
\end{equation}
Since $\deg C_i\le 1$, we have $\deg D\le 3$. Write
\begin{gather*}
D=\sum_{i=0}^3D_iY^i,\cr
C_1^2C_4+C_1C_2C_3+C_0C_3^2=\sum_{i=0}^3E_iY^i,\cr
(C_1^2C_4+C_1C_2C_3+C_0C_3^2)^2+(C_1C_3+C_2^2)^3=\sum_{i=0}^6F_iY^i,
\end{gather*}
where
\begin{align}
\label{E0}
E_0=\,& k (a^5+a^4+a^3 b_1^2+a^3+a^2 b_1^2+a^2)+k^3
(a b_1^4+b_1^4)\\ 
&+k^2 (a b_1^4+a
b_1^2+b_1^5+b_1^4+b_1^3+b_1^2), \cr 
\label{E1}
E_1=\,& a^5+a^4+a^3 b_1^2+a^3+a^2 b_1^2+a^2+k (a^4
b_1+a^2 b_1^3+a^2 b_1)\\
&+k^2 (a
b_1^4+b_1^5+b_1^4+b_1^3)+ k^3b_1^5, \cr
\label{E2}
E_2=\,& a^5+a^4 b_1+a^4+a^3+a^2 b_1+a^2+k (a
b_1^2+b_1^2)+k^2 (a
b_1^4+b_1^5+b_1^4), \\ 
\label{E3}
E_3=\,& a^4 b_1+a^2 b_1+a b_1^2+b_1^2+kb_1^3  +k^2b_1^5. 
\end{align}
\begin{align}
\label{F0}
F_0=&\,k^2
(a^{10}+a^8+a^6 b_1^4+a^6+a^4 b_1^4+a^4)+ k^3
(a^{12}+a^{10} b_1+a^{10}+a^9 b_1^3+a^9
b_1\\
&+a^8 b_1^2+a^8+a^6 b_1^6+a^6 b_1^3+a^6
b_1+a^6+a^5 b_1^5+a^5 b_1+a^4 b_1^7+a^4
b_1^6+a^4 b_1^3\cr
&+a^4 b_1^2+a^3 b_1^9+a^3
b_1^7+a^3 b_1^5+a^3 b_1^3)+  k^4
(a^8 b_1^3+a^8 b_1^2+a^4 b_1^5+a^4
b_1^4+a^4 b_1^3\cr
&+a^4 b_1^2+a^2 b_1^9+a^2
b_1^5+b_1^{10}+b_1^8+b_1^6+b_1^4)+k^5
(a^8 b_1^4+a^2 b_1^{10}+a^2 b_1^7+a^2
b_1^5\cr
&+a^2 b_1^4+a b_1^9+a b_1^5)+ k^6 (a^2 b_1^8+b_1^9+b_1^7+b_1^6)+ k^7
(a^4 b_1^8+a^2 b_1^9+a^2 b_1^8+a
b_1^{11}\cr
&+a b_1^9+b_1^{10}+b_1^8) + k^8(b_1^{11}+b_1^{10}) +k^9 b_1^{12}, \cr 
\label{F1}
F_1=&\,k^2 (a^{12}+a^{11}
b_1+a^{10} b_1+a^{10}+a^8 b_1^2+a^8+a^7
b_1^3+a^7 b_1+a^6 b_1^6+a^6 b_1^3+a^6
b_1\\
&+a^6+a^5 b_1^7+a^5 b_1^3+a^4 b_1^7+a^4
b_1^6+a^4 b_1^3+a^4 b_1^2)+ k^3 (a^9
b_1^3+a^8 b_1^3+a^8 b_1^2\cr
&+a^5 b_1^5+a^5
b_1^3+a^4 b_1^5+a^4 b_1^4+a^4 b_1^3+a^4
b_1^2+a^3 b_1^9+a^3 b_1^5+a^2 b_1^9+a^2
b_1^8+a^2 b_1^5\cr
&+a^2 b_1^4)+ k^4 (a^8 b_1^4+a^3 b_1^7+a^3 b_1^5+a^2
b_1^{10}+a^2 b_1^7+a^2 b_1^5+a^2 b_1^4)+ k^5 (a b_1^9+a
b_1^7\cr
&+b_1^9+b_1^8+b_1^7+b_1^6)+ k^6
(a^4 b_1^8+a^3 b_1^9+a^2 b_1^9+a^2
b_1^8+b_1^{10}+b_1^8)\cr
&+ k^7 (a
b_1^{11}+b_1^{11}+b_1^{10})+ k^8 b_1^{
12}, \cr 
\label{F2}
F_2=&\, a^{10}+a^8+a^6 b_1^4+a^6+a^4 b_1^4+a^4+k (a^{12}+a^{10}
b_1^2+a^{10} b_1+a^{10}+a^9 b_1^3+a^9 b_1\\
&+a^8
b_1^3+a^8+a^7 b_1^5+a^7 b_1^3+a^6 b_1^3+a^6
b_1^2+a^6 b_1+a^6+a^5 b_1^5+a^5 b_1)+k^2
(a^{12}\cr
&+a^{10} b_1+a^{10}+a^9 b_1+a^8
b_1^3+a^8 b_1^2+a^8+a^6 b_1^6+a^6 b_1^5+a^6
b_1^4+a^6 b_1^3+a^6 b_1\cr
&+a^6+a^5 b_1^3+a^5
b_1+a^4 b_1^7+a^4 b_1^5+a^4 b_1^4+a^4
b_1^2+a^3 b_1^7+a^3 b_1^3)+ k^3
(a^8 b_1^4\cr
&+a^8 b_1^3+a^8 b_1^2+a^4
b_1^7+a^4 b_1^6+a^4 b_1^5+a^4 b_1^4+a^4
b_1^3+a^4 b_1^2+a^3 b_1^9+a^3 b_1^7+a^2
b_1^9\cr
&+a^2 b_1^8+a^2 b_1^7+a^2 b_1^6+a
b_1^9+a b_1^5)+ k^4
(a^8 b_1^4+a^2 b_1^{10}+a^2 b_1^9+a^2
b_1^7+a^2 b_1^5\cr
&+a^2 b_1^4+a b_1^7+a
b_1^5+b_1^{10}+b_1^9+b_1^7)+k^5
(a^4 b_1^8+a^2 b_1^{10}+a^2
b_1^9+a^2 b_1^8+a b_1^{11}\cr
&+a
b_1^9+b_1^{10}+b_1^9+b_1^7+b_1^6)+ k^6
(a^4 b_1^8+a^2 b_1^9+a^2 b_1^8+a
b_1^9+b_1^{11}+b_1^{10}+b_1^8)\cr
&+ k^7(b_1^{12}+b_1^{11}+b_1^{10}) + k^8 b_1^{12}, \cr 
\label{F3}
F_3=&\, a^{12}+a^{11} b_1+a^{10} b_1^2+a^{10} b_1+a^{10}+a^9
b_1^3+a^8 b_1^3+a^8+a^7 b_1^3+a^7 b_1+a^6
b_1^3\\
&+a^6 b_1^2+a^6 b_1+a^6+ k (a^9 b_1^3+a^8
b_1^3+a^8 b_1^2+a^7 b_1^5+a^6 b_1^5+a^6
b_1^4+a^5 b_1^5+a^5 b_1^3\cr
&+a^4 b_1^5+a^4
b_1^4+a^4 b_1^3+a^4 b_1^2)+ k^2
(a^8 b_1^4+a^5 b_1^7+a^4 b_1^7+a^4 b_1^6+a^3
b_1^7+a^3 b_1^5\cr
&+a^2 b_1^7+a^2 b_1^6+a^2
b_1^5+a^2 b_1^4)+ k^3 (a^3
b_1^9+a^2 b_1^9+a^2 b_1^8+a b_1^9+a
b_1^7+b_1^9+b_1^8\cr
&+b_1^7+b_1^6)+k^4
(a^4 b_1^8+a^3 b_1^9+a^2 b_1^{10}+a^2
b_1^9+a^2 b_1^8+b_1^{10}+b_1^8)\cr
&+k^5 (a
b_1^{11}+b_1^{11}+b_1^{10})+  k^6 b_1^{12}, \cr 
\label{F4}
F_4=&\, a^{12}+a^{10} b_1^2+a^{10} b_1+a^9 b_1+a^8 b_1^3+a^8
b_1^2+a^7 b_1^3+a^6 b_1^3+a^6 b_1^2+a^6
b_1+a^5 b_1^3\\
&+a^5 b_1+a^4 b_1^2+a^4+ k
(a^{12}+a^{10} b_1+a^{10}+a^9 b_1^3+a^9 b_1+a^8
b_1^3+a^8+a^6 b_1^5\cr
&+a^6 b_1^4+a^6 b_1^3+a^6
b_1+a^6+a^5 b_1^5+a^5 b_1+a^4 b_1^5+a^4
b_1^4+a^3 b_1^5+a^3 b_1^3)\cr
&+ k^2
(a^8 b_1^4+a^8 b_1^3+a^8 b_1^2+a^4 b_1^7+a^4
b_1^6+a^4 b_1^5+a^4 b_1^4+a^4 b_1^3+a^4
b_1^2+a^3 b_1^7\cr
&+a^2 b_1^7+a^2 b_1^6+a^2
b_1^4+a b_1^7+a b_1^5+b_1^4)+ k^3
(a^8 b_1^4+a^2 b_1^9+a^2 b_1^8+a^2
b_1^7+a^2 b_1^5\cr
&+a^2 b_1^4+a b_1^9+a
b_1^5+b_1^9+b_1^8+b_1^7+b_1^6)+ k^4
(a^4 b_1^8+a^2 b_1^{10}+a^2 b_1^9+a
b_1^9+b_1^9\cr
&+b_1^8+b_1^7+b_1^6)+ k^5
(a^4 b_1^8+a^2 b_1^9+a^2 b_1^8+a
b_1^{11}+a b_1^9+b_1^{11}+b_1^8)\cr
&+ k^6(b_1^{12}+b_1^{11}+b_1^{10}) +k^7 b_1^{12}, \cr 
\label{F5}
F_5=&\, a^{12}+a^{11} b_1+a^{10} b_1+a^{10}+a^8 b_1^2+a^8+a^7
b_1^3+a^7 b_1+a^6 b_1^3+a^6 b_1+a^6+a^5
b_1^3\\
&+a^4 b_1^3+a^4 b_1^2+ k (a^9 b_1^3+a^8
b_1^3+a^8 b_1^2+a^5 b_1^5+a^5 b_1^3+a^4
b_1^5+a^4 b_1^4+a^4 b_1^3\cr
&+a^4 b_1^2+a^3
b_1^5+a^2 b_1^5+a^2 b_1^4)+ k^2
(a^8 b_1^4+a^3 b_1^7+a^3 b_1^5+a^2
b_1^7+a^2 b_1^5+a^2 b_1^4)\cr
&+k^3 (a b_1^9+a
b_1^7+b_1^9+b_1^8+b_1^7+b_1^6)+ k^4
(a^4 b_1^8+a^3 b_1^9+a^2 b_1^9+a^2
b_1^8+b_1^{10}+b_1^8)\cr
&+k^5 (a
b_1^{11}+b_1^{11}+b_1^{10})+  k^6 b_1^{
12}, \cr 
\label{F6}
F_6=&\, a^{12}+a^{10} b_1+a^{10}+a^9 b_1+a^8+a^6
b_1^3+a^6 b_1+a^6+a^5 b_1^3+a^5 b_1+a^4
b_1^3+a^3 b_1^3\\
& +b_1^4 +a^2 b_1^4+k (a^8 b_1^3+a^8
b_1^2+a^4 b_1^5+a^4 b_1^4+a^4 b_1^3+a^4
b_1^2+a^2 b_1^5+a^2 b_1^4)\cr
&+ k^2 (a^8 b_1^4+a^2
b_1^7+a^2 b_1^5+a^2 b_1^4+a b_1^7+a
b_1^5+b_1^6)+  k^3(b_1^9+b_1^8+b_1^7+b_1^6)\cr
& +k^4 (a^4 b_1^8+a^2 b_1^9+a^2
b_1^8+a b_1^9+b_1^8) + k^5 (b_1^{11}+b_1^{10})
+k^6 b_1^{12}. \nonumber
\end{align}
Then \eqref{2.9} becomes
\[
\Bigl(\sum_{i=0}^3D_iY^i\Bigr)\Bigl(\sum_{j=0}^3(D_j+E_j)Y^j\Bigr)=\sum_{i=0}^6F_iY^i.
\]
Comparing the coefficients in the above gives
\begin{align}
\label{2.10}
&D_3^2+D_3E_3=F_6,\\
\label{2.11}
&D_3E_2+D_2E_3=F_5,\\
\label{2.12}
&D_2^2+D_3E_1+D_2E_2+D_1E_3=F_4,\\
\label{2.13}
&D_3E_0+D_2E_1+D_1E_2+D_0E_3=F_3,\\
\label{2.14}
&D_1^2+D_2E_0+D_1E_1+D_0E_2=F_2,\\
\label{2.15}
&D_1E_0+D_0E_1=F_1,\\
\label{2.16}
&D_0^2+D_0E_0=F_0.
\end{align}
The above is a system of seven polynomial equations in seven unknowns $D_0,D_1,D_2$, $D_3,a,b_1,k$ in $\f_q$. In what follows, we show that the above system does not allow a solution with $a,b_1,k\in\f_q^*$ such that $\text{Tr}_{q/2}(k)=1$ and $1+aX+bX^3$ has no root in $\mu_{q+1}$. This proves Claim~\ref{C3.1} by contradiction.

Before proceeding, we mention that $E_2$ and $E_3$ are not both $0$, a claim that was used earlier in this subsection to transform \eqref{2.6} to \eqref{2.7}. We have $\text{Res}(E_2,E_3;a)=b_1^{17}k^2\ne 0$, hence the claim.

\subsection{Nonsolvability of the system}

$1^\circ$ We temporarily assume that $E_3\ne 0$. By \eqref{2.11} -- \eqref{2.13},
\begin{align*}
D_2\,&=(D_3E_2+F_5)/E_3,\cr
D_1\,&=(D_2^2+D_3E_1+D_2E_2+F_4)/E_3,\cr
D_0\,&=(D_3E_0+D_2E_1+D_1E_2+F_3)/E_3.
\end{align*}
By these substitutions, the remaining equations in \eqref{2.10} -- \eqref{2.16} become 
\begin{align}
\label{2.17}
&D_3^2+D_3E_3+F_6=0,\\
\label{2.18}
& D_3^4 E_2^4+D_3^2 E_1^2 E_3^4+D_3^2
E_1 E_2^2 E_3^3+D_3 E_1^2 E_3^5+D_3
E_1 E_2^2 E_3^4+D_3 E_2^4
E_3^3+E_0 E_3^5 F_5\\
&+E_1 E_3^5
F_4+E_1 E_3^3 F_5^2+E_2^3 E_3^3
F_5+E_2^2 E_3^4 F_4+E_2 E_3^5
F_3+E_3^6 F_2+E_3^4 F_4^2+F_5^4 =0,\cr
\label{2.19}
& D_3^2 E_0 E_2^2 E_3+D_3^2 E_1
E_2^3+D_3 E_0 E_2^2 E_3^2+D_3
E_1 E_2^3 E_3+E_0 E_2 E_3^2
F_5+E_0 E_3^3 F_4\\
&+E_0 E_3
F_5^2+E_1^2 E_3^2 F_5+E_1 E_2^2
E_3 F_5+E_1 E_2 E_3^2 F_4+E_1
E_2 F_5^2+E_1 E_3^3 F_3+E_3^4 F_1 =0,\cr
\label{2.20}
& D_3^4 E_2^6+D_3^2 E_0^2 E_3^6+D_3^2
E_0 E_2^3 E_3^4+D_3^2 E_2^6
E_3^2+D_3 E_0^2 E_3^7+D_3 E_0
E_2^3 E_3^5+E_0 E_1 E_3^6 F_5\\
&+E_0
E_2^2 E_3^5 F_5+E_0 E_2 E_3^6
F_4+E_0 E_2 E_3^4 F_5^2+E_0
E_3^7 F_3+E_1^2 E_3^4 F_5^2+E_2^4
E_3^2 F_5^2\cr
&+E_2^2 E_3^4 F_4^2+E_2^2
F_5^4+E_3^8 F_0+E_3^6 F_3^2 =0. \nonumber
\end{align}
In fact, \eqref{2.17} --  \eqref{2.20} are linear combinations of \eqref{2.10} -- \eqref{2.16} with coefficients in $\f_2[D_0,\dots,D_3,E_0,\dots,E_3,F_0,\dots,F_6]$. Hence they hold even if $E_3=0$. Thus we remove the assumption that $E_3\ne 0$. Subtracting \eqref{2.18} -- \eqref{2.20} by suitable multiples of \eqref{2.17}  allows us to eliminate $D_3$, resulting in three equations in $E_i$'s and $F_j$'s:
\begin{align}
\label{2.21}
H_1:=\,& E_0 E_3^5 F_5+E_1^2 E_3^4
F_6+E_1 E_2^2 E_3^3 F_6+E_1
E_3^5 F_4+E_1 E_3^3 F_5^2+E_2^4
E_3^2 F_6+E_2^4 F_6^2\\
&+E_2^3 E_3^3
F_5+E_2^2 E_3^4 F_4+E_2 E_3^5
F_3+E_3^6 F_2+E_3^4 F_4^2+F_5^4 =0, \nonumber \\
\label{2.22}
H_2:=\,& E_0 E_2^2 E_3 F_6+E_0 E_2
E_3^2 F_5+E_0 E_3^3 F_4+E_0
E_3 F_5^2+E_1^2 E_3^2 F_5+E_1
E_2^3 F_6\\
&+E_1 E_2^2 E_3
F_5+E_1 E_2 E_3^2 F_4+E_1 E_2
F_5^2+E_1 E_3^3 F_3+E_3^4 F_1 =0, \nonumber \\
\label{2.23}
H_3:=\,& E_0^2 E_3^6 F_6+E_0 E_1 E_3^6
F_5+E_0 E_2^3 E_3^4 F_6+E_0 E_2^2
E_3^5 F_5+E_0 E_2 E_3^6 F_4\\
&+E_0
E_2 E_3^4 F_5^2+E_0 E_3^7 F_3+E_1^2
E_3^4 F_5^2+E_2^6 F_6^2+E_2^4 E_3^2
F_5^2+E_2^2 E_3^4 F_4^2\cr
&+E_2^2
F_5^4+E_3^8 F_0+E_3^6 F_3^2 =0. \nonumber
\end{align}
$H_1,H_2,H_3$ are polynomials in $a,b_1,k$:
\begin{align}
\label{2.24}
H_1\,&= C_4^4h_1,\\
\label{2.25}  
H_2\,&= C_4^5(a^2+b_1^2k)^2h_2,\\
\label{2.26}
H_3\,&= C_4^4h_3,
\end{align}
where $h_1,h_2,h_3$ are given in Appendix \eqref{A1} -- \eqref{A3}. We note that $\deg_kh_1=20$, $\deg_kh_2=10$, and $\deg_kh_3=24$.

We claim that $a^2+b_1^2k\ne 0$. Assume to the contrary that $k=a^2/b_1^2$. Then 
\begin{align*}
h_1\,&=a^7b_1^{10}(a^3+b_1+a^2b_1),\cr
h_3\,&=a^9b_1^{12}(a^5+b_1+a^4b_1+b_1^3+a^4b_1^3),
\end{align*}
and
\[
\text{Res}(a^3+b_1+a^2b_1,\, a^5+b_1+a^4b_1+b_1^3+a^4b_1^3;\, b_1)=a^3(1+a)^4(1+a+a^3)^2.
\]
Thus $(1+a)(1+a+a^3)=0$. If $1+a=0$, then $a^3+b_1+a^2b_1=1\ne 0$, which is a contradiction. If $1+a+a^3=0$, then from $a^3+b_1+a^2b_1=0$ we have $b_1=a^3/(1+a)^2=a^{-3}$, hence $k=a^2/b_1^2=a^8$. Note that $a\in\f_{2^3}\subset\f_q$ and $\text{Tr}_{2^3/2}(a)=0$. Thus $\text{Tr}_{q/2}(k)=\text{Tr}_{q/2}(a)=0$, which is a contradiction. Thus the claim is proved.

Now by \eqref{2.24} -- \eqref{2.26}, we have $h_1=h_2=h_3=0$. Using $h_2$ to reduce the degree of $k$ in $h_1$ and $h_3$, we find that $h_1\equiv b_1^8h_1'\pmod{h_2}$ and $h_3\equiv b_1^{10}h_3'\pmod{h_2}$, where $h_1'$ and $h_3'$, given in \eqref{A4} and \eqref{A5}, are of degree $9$ in $k$. To recap, we now have three polynomial equations in $a,b_1,k$:
\begin{equation}\label{eq3.36}
h_1'=0,\quad h_2=0,\quad h_3'=0.
\end{equation}

\medskip
$2^\circ$ We have
\begin{align}
\label{2.26.1}
\text{Res}(h_1',h_2;a)\,&= b_1^{208}k^{37}(1+b_1+b_1^2k)^2(1+b_1^2k)^2S_1^6S_2,\\
\label{2.26.2}
\text{Res}(h_2,h_3';a)\,&= b_1^{272}k^{69}(1+b_1+b_1^2k)^2S_1^8S_3,\\
\label{2.27}
\text{Res}(h_1',h_2;b_1)\,&= a^{108}(1+a)^{108}k^{131}T_1^6T_2, \\
\label{2.28}
\text{Res}(h_2,h_3';b_1)\,&= a^{126}(1+a)^{152}k^{158}T_1^8T_3,
\end{align}
where $S_1,S_2,S_3,T_1,T_2,T_3$ are given in \eqref{A6} -- \eqref{A11}. In \eqref{2.26.1} and \eqref{2.26.2}, $1+b_1+b_1^2k\ne0$ since $\text{Tr}_{q/2}(k)=1$.

\medskip
{\bf Case 1.} Assume $a=1$. We have $h_1'|_{a=1}=b_1^9k^3h_1''$ and $h_2|_{a=1}=b_1^9k^3h_3''$, where
\begin{align*}
h_1''=\,& b_1^{10}
k^6+ b_1^{11} k^5+b_1^9 k^5+b_1^8 k^5+ b_1^{10} k^4+ b_1^9
k^4+b_1^9 k^3+b_1^7 k^2\cr
&+b_1^6
k^2+b_1^2 k^2+b_1 k+k+ b_1^4+b_1^3+b_1^2,\cr
h_2''=\,& b_1^{11} k^7+b_1^9 k^6+b_1^8 k^5+b_1^8
k^4+b_1^7 k^3+b_1^6 k^3+b_1^3 k^3+ b_1^5 k^2\cr
&+b_1^4
k^2+b_1^2 k+b_1 k^2+k +b_1^2.
\end{align*}
It turns out that $\text{Res}(h_1'',h_2'';b_1)=k^{66}\ne 0$, which is a contradiction.

\medskip
{\bf Case 2.} Assume that $a\ne 1$ but $T_1=0$, i.e., $k=a^{-2}+a^{-1}+1+a^3+a^4+a^6=(1+a)^5(1+a^2+a^3)/a^2$. We have
\begin{align}
\label{h1T1}
h_1'|_{T_1=0}\,&=a^{-16}(1+a)^4(a+b_1+ab_1+a^2b_1+a^3b_1)h_1^*,\\
\label{h2T1}
h_2|_{T_1=0}\,&=a^{-20}(1+a)^4(a+b_1+ab_1+a^2b_1+a^3b_1)h_2^*,\\
\label{h3T1}
h_3'|_{T_1=0}\,&=a^{-17}(1+a)^4(a+b_1+ab_1+a^2b_1+a^3b_1)h_3^*.\\
\end{align}
where $h_1^*,h_2^*,h_3^*$ are given in \eqref{A12} -- \eqref{A14}.

First assume that $a+b_1+ab_1+a^2b_1+a^3b_1\ne 0$. Then $h_1^*=h_2^*=h_3^*=0$. However, we find that as a polynomial in $b_1$, 
\[
\text{gcd}(\text{Res}(h_1^*,h_2^*;a),\text{Res}(h_2^*,h_3^*;a))=b_1^{1720},
\]
which dose not have a nonzero root. This is a contradiction.

Now assume that $a+b_1+ab_1+a^2b_1+a^3b_1= 0$, i.e., $b_1=a/(1+a)^3$. Making substitutions $k=(1+a)^5(1+a^2+a^3)/a^2$ and $b_1=a/(1+a)^3$ gives $E_3=0$.

\medskip 
{\bf Case 3.} Assume that $(1+a)T_1\ne0$. Then by \eqref{2.27} and \eqref{2.28}, $T_2=T_3=0$. We find that 
\begin{align}
\label{2.29}
\text{Res}(T_2,T_3;a)\,&=T_4^2T_5^2T_6^2T_7^2,\\
\label{2.30}
\text{Res}(T_2,T_3;k)\,&=(1+a)^{44}T_8^2T_9^2T_{10}^2T_{11}^2,
\end{align}
where
\begin{align*}
T_4=\,&k^7+k^6+k^5+k^3+k^2+k+1, \cr
T_5=\,&k^{13}+k^{12}+k^{11}+k^{10}+k^8+k^7+k^6+k^5+k^4+k^3+1, \cr
T_6=\,&k^{21}+k^{20}+k^{19}+k^{18}+k^{16}+k^{14}+k^{12}+k^{11}+k^{10}+k^8+k^7+k^6+k^5\cr
&+k^4+k^3+k^2+1, \cr
T_7=\,&k^{49}+k^{46}+k^{44}+k^{40}+k^{39}+k^{38}+k^{37}+k^{36}+k^{35}+k^{30}+k^{24}+k^{
22}+k^{15}\cr
&+k^{14}+k^{12}+k^{10}+k^8+k^3+k^2+k+1, \cr
T_8=\,&a^7+a^3+1, \cr
T_9=\,&a^{13}+a^{11}+a^{10}+a^5+a^3+a+1, \cr
T_{10}=\,&a^{21}+a^{20}+a^{18}+a^{14}+a^{11}+a^8+a^7+a^4+1, \cr
T_{11}=\,&a^{49}+a^{48}+a^{47}+a^{44}+a^{43}+a^{42}+a^{40}+a^{38}+a^{37}+a^{36}+a^{33}+a^{
31}+a^{26}\cr
&+a^{25}+a^{22}+a^{20}+a^{19}+a^{18}+a^{17}+a^{16}+a^{15}+a^8+a^6+a^5
+1. 
\end{align*}
In the above $T_7\ne 0$ since $\text{Tr}_{q/2}(k)=1$. 

\medskip
{\bf Case 3.1.} Assume in \eqref{2.26.1} that $(1+b_1^2k)S_1\ne 0$. Then by \eqref{2.26.1} and \eqref{2.26.2}, $S_2=S_3=0$. We find that, as a polynomial in $k$,
\[
\text{gcd}(\text{Res}(S_2,S_3;b_1),\text{Res}(T_2,T_3;a))=1,
\]
which is a contradiction. (Recall that $\text{Res}(T_2,T_3;a)$ is given in \eqref{2.29}.)

\medskip
{\bf Case 3.2.} Assume in \eqref{2.26.1} that $(1+b_1^2k)S_1= 0$.

First assume that $1+b_1^2k=0$. We have $h_1'|_{k=b_1^{-2}}=ah_1^\dagger$, $h_2|_{k=b_1^{-2}}=ah_2^\dagger$, $h_3'|_{k=b_1^{-2}}=ah_3^\dagger$, where
\begin{align*}
h_1^\dagger=\,& a^{21}+a^{19} b_1^4+a^{19} b_1^2+a^{17}+a^{16} b_1^5+a^{16}
b_1^3+a^{14} b_1^5+a^{14} b_1^3+a^{13} b_1^6+a^{13}
b_1^4\cr
&+a^{13}+a^{12} b_1^5+a^{11} b_1^6+a^{11}
b_1^2+a^{10} b_1^5+a^9 b_1^6+a^9 b_1^4+a^9+a^8
b_1^7+a^8 b_1^3\cr
&+a^6 b_1^7+a^6 b_1^5+a^6
b_1^3+a^4 b_1^5+a^2 b_1^5+a b_1^4+b_1^7,\cr
h_2^\dagger=\,& a^{19}+a^{17}+a^{16} b_1+a^{14} b_1^3+a^{12} b_1^3+a^{11}
b_1^4+a^{11}+a^{10} b_1^3+a^9 b_1^4+a^9+a^8
b_1^5\cr
&+a^8 b_1^3+a^8 b_1+a^7 b_1^4+a^6
b_1^5+a^5 b_1^4+a^3 b_1^4+a b_1^4+b_1^5,\cr
h_3^\dagger=\,& a^{23} b_1^2+a^{23}+a^{22} b_1^3+a^{21} b_1^4+a^{21}
b_1^2+a^{21}+a^{20} b_1^5+a^{20} b_1+a^{19}
b_1^4+a^{19}\cr
&+a^{18} b_1^5+a^{18} b_1^3+a^{17}
b_1^6+a^{17} b_1^4+a^{17} b_1^2+a^{17}+a^{16}
b_1^7+a^{16} b_1^5+a^{16} b_1\cr
&+a^{15}
b_1^8+a^{15} b_1^6+a^{15} b_1^4+a^{15}
b_1^2+a^{15}+a^{14} b_1^5+a^{13} b_1^6+a^{13}
b_1^4+a^{13} b_1^2\cr
&+a^{13}+a^{12} b_1^9+a^{12}
b_1^7+a^{12} b_1^5+a^{12} b_1^3+a^{12}
b_1+a^{11} b_1^8+a^{11} b_1^4+a^{11}\cr
&+a^{10}
b_1^7+a^{10} b_1^5+a^9 b_1^{10}+a^9 b_1^8+a^9
b_1^4+a^9 b_1^2+a^9+a^8 b_1^9+a^8 b_1^7+a^8
b_1^5\cr
&+a^8 b_1^3+a^8 b_1+a^7 b_1^{10}+a^7
b_1^4+a^5 b_1^6+a^4 b_1^9+a^4 b_1^7+a^4
b_1^5+a^3 b_1^6+a^3 b_1^4\cr
&+a b_1^{10}+a
b_1^8+a b_1^6+a b_1^4+b_1^9+b_1^5.
\end{align*}
Moreover,
\begin{align}
\label{2.32}
\text{Res}(h_1^\dagger,h_2^\dagger;a)\,&=b_1^{119}(1+b_1+b_1^2)^{10}(1+b_1+b_1^4)^2,\\
\label{2.33}
\text{Res}(h_2^\dagger,h_3^\dagger;a)\,&=b_1^{135}(1+b_1+b_1^2)^{8}(1+b_1+b_1^2+b_1^4+b_1^6)^2,\\
\label{2.34}
\text{Res}(h_1^\dagger,h_2^\dagger;b_1)\,&=a^{65}(1+a)^{54}(1+a+a^2)^{10}(1+a+a^4)^2(1+a^3+a^4)^8.
\end{align}
By \eqref{2.32} and \eqref{2.33}, we have $1+b_1+b_1^2=0$. Thus $k=b_1^{-2}\in\f_{2^2}$. Since $\text{Tr}_{q/2}(k)=1$, $[\f_q:\f_{2^2}]$ is odd. Then by \eqref{2.34}, we have $1+a+a^2=0$. Therefore, $a=b_1$ or $b_1^{-1}$. However,
\[
h_2^\dagger|_{a=b_1^{-1}}=b_1^{-19}(1+b_1+b_1^3)^2(1+b_1^2+b_1^3)^2(1+b_1^3+b_1^6)^2\ne 0.
\]
So we must have $a=b_1$. Recall that $\text{Tr}_{q^2/q}(b)=b_1$ and $\text{N}_{q^2/q}(b)=b_1^2k=1$. Hence $b^2+b_1b+1=0$. Then $b^{-1}\in\mu_{q+1}$ is a root of $1+aX+bX^3$, which is a contradiction.

Now assume that $S_1=0$. Recall that in \eqref{2.29}, $T_7\ne 0$. Hence $T_4T_5T_6=0$. Then $[\f_2(k):\f_2]$ is odd. Since $\text{Tr}_{q/2}(k)=1$, $[\f_q:\f_2]$ is also odd.

\medskip
{\bf Case 3.2.1.} Assume that $T_4=0$. We have 
\[
\text{Res}(S_1,T_4;k)=T_{12}T_{13}T_{14},
\]
where
\begin{align*}
T_{12}=\,& b_1^7+b_1^6+b_1^5+b_1^4+b_1^3+b_1^2+1,\cr
T_{13}=\,& b_1^{14}+b_1^{13}+b_1^{12}+b_1^{11}+b_1^{10}+b_1^8+1,\cr
T_{14}=\,& b_1^{35}+b_1^{33}+b_1^{31}+b_1^{29}+b_1^{28}+b_1^{27}+b_1^{26}+b_1^{25}+b_1^{24}+b_1^{23}+b_1
^{21}+b_1^{20}\cr
&+b_1^{17}+b_1^{14}+b_1^{13}+b_1^{
12}+b_1^9+b_1^8+b_1^7+b_1+1.
\end{align*}
Since $[\f_q:\f_2]$ is odd, $T_{13}\ne0$. Thus $T_{12}T_{14}=0$. 

First assume that $T_{12}=0$, i.e., $T_{12}$ is the minimal polynomial of $b_1$ over $\f_2$. Treat $S_1$ and $T_4$ as polynomials in $k$ over $\f_2(b_1)$, where $\deg_k S_1=3$ and $\deg_kT_4=7$. Applying the Euclidean algorithm to $S_1$ and $T_4$ gives some $r(k)\in\f_2(b_1)[k]$ such that $r(k)\in\f_2(b_1)S_1+\f_2(b_1)T_4$ and $\deg r(k)=1$; solving $r(k)=0$ gives
\[
k=b_1+b_1^3+b_1^4.
\]
Now $h_1',h_2,h_3'$ can be expressed as polynomials in $a$ and $b_1$ which are further reduced modulo $T_{12}$. We find that, as a polynomial in $a$,
\[
\text{gcd}(\text{Res}(h_1',h_2;b_1),\text{Res}(T_2,T_3;k))=1+a.
\]
(Recall that $\text{Res}(T_2,T_3;k)$ is given in \eqref{2.30}.) Since $a\ne 1$ (by Case 1), we have a contradiction. 

Next, assume that $T_{14}=0$, i.e., $T_{14}$ is the minimal polynomial of $b_1$ over $\f_2$. By the same computation, we have 
\[
k=b_1^2+ b_1^4 + b_1^6 +b_1^7+b_1^{10}+b_1^{11}+b_1^{18}+b_1^{21}+b_1^{22}+b_1^{24}+b_1^{25}+b_1^{26}+b_1^{27}+b_1^{30}.
\]
and, as a polynomial in $a$,
\[
\text{gcd}(\text{Res}(h_1',h_2;b_1),\text{Res}(T_2,T_3;k))=(1+a)^3,
\]
which is a contradiction.

\medskip
{\bf Case 3.2.2.} Assume that $T_5=0$. We have 
\[
\text{Res}(S_1,T_5;k)=T_{15}T_{16}T_{17},
\]
where
\begin{align*}
T_{15}=\,& b_1^{13}+b_1^{12}+b_1^6+b_1^5+b_1^4+b_1^3+b_1^2+b_1+1,\cr
T_{16}=\,& b_1^{13}+b_1^{12}+b_1^{10}+b_1^9+b_1^5+b_1^3+
b_1^2+b_1+1,\cr
T_{17}=\,& b_1^{78}+b_1^{76}+b_1^{75}+b_1^{74}+b_1^{72}+b_1^{69}+b_1^{68}+b_1^{67}+b_1^{66}+b_1^{65}+b_1
^{64}+b_1^{63}+b_1^{62}\cr
&+b_1^{61}+b_1^{60}+b_1^{
59}+b_1^{58}+b_1^{55}+b_1^{54}+b_1^{53}+b_1^{52}
+b_1^{51}+b_1^{50}+b_1^{49}+b_1^{46}\cr
&+b_1^{44}+
b_1^{41}+b_1^{39}+b_1^{38}+b_1^{36}+b_1^{35}+
b_1^{32}+b_1^{31}+b_1^{30}+b_1^{29}+b_1^{28}+
b_1^{26}+b_1^{25}\cr
&+b_1^{22}+b_1^{20}+b_1^{19}+
b_1^{18}+b_1^{15}+b_1^{14}+b_1^{13}+b_1^9+
b_1^6+b_1^5+b_1^2+b_1+1.
\end{align*}
Since $[\f_q:\f_2]$ is odd, $T_{17}\ne0$. Thus $T_{15}T_{16}=0$. We now go through the same computation as in Case~3.2.1 (with $T_4$ replaced by $T_5$):

When $T_{15}=0$, we have
\[
k=b_1^3+b_1^6+b_1^8+b_1^9,
\]
and, as a polynomial in $a$,
\[
\text{gcd}(\text{Res}(h_1',h_2;b_1),\text{Res}(T_2,T_3;k))=1,
\]
which is a contradiction.

When $T_{16}=0$, we have
\[
k=b_1^2+b_1^3+b_1^4+b_1^6+b_1^{11},
\]
and, as a polynomial in $a$,
\[
\text{gcd}(\text{Res}(h_1',h_2;b_1),\text{Res}(T_2,T_3;k))=1+a,
\]
which is a contradiction. 

\medskip
{\bf Case 3.2.3.} Assume that $T_6=0$. We have 
\[
\text{Res}(S_1,T_6;k)=T_{18}T_{19}T_{20},
\]
where
\begin{align*}
T_{18}=\,&b_1^{21}+b_1^{19}+b_1^{18}+b_1^{17}+b_1^{15}+b_1^{13}+b_1^{12}+b_1^{10}+b_1^9+b_1^8+b_1^7+b_1^6+1, \cr
T_{19}=\,& b_1^{21}+b_1^{20}+b_1^{19}+b_1^{18}+b_1^{17}+b_1^{15}+b_1^{14}+b_1^8+b_1^7+b_1^6+b_1^4+b_1^3+1,\cr
T_{20}=\,& b_1^{126}+b_1^{125}+b_1^{124}+b_1^{122}+b_1^{120}+b_1^{119}+b_1^{117}+b_1^{116}+b_1^{114}+b_1^{112}
+b_1^{109}\cr
&+b_1^{105}+b_1^{104}+b_1^{103}+b_1^{
100}+b_1^{99}+b_1^{95}+b_1^{93}+b_1^{90}+b_1^{88
}+b_1^{87}+b_1^{86}+b_1^{85}\cr
&+b_1^{84}+b_1^{80}+b_1^{79}+b_1^{72}+b_1^{71}+b_1^{70}+b_1^{69}+b_1^{68}+b_1^{67}+b_1^{64}+b_1^{63}+b_1^{60}+b_1^{59}\cr
&+b_1^{56}+b_1^{52}+b_1^{50}+b_1^{47}+b_1^{45}+b_1^{43}+b_1^{39}+b_1^{37}+b_1^{34}+b_1^{33}+b_1^{32}+b_1^{27}+b_1^{26}\cr
&+b_1^{20}+b_1^{19}+b_1^{15}+b_1^{13}+b_1^{12}+b_1^6+b_1^5+b_1^4+b_1+1.
\end{align*}
Since $[\f_q:\f_2]$ is odd, $T_{20}\ne0$. Thus $T_{18}T_{19}=0$. Again, we go through the same computation as in Case~3.2.1 (with $T_4$ replaced by $T_6$):

When $T_{18}=0$, we have
\[
k=1+b_1+b_1^2+b_1^3+b_1^5+b_1^6+b_1^8+b_1^{11}+b_1^{12}+b_1^{14}+b_1^{15}+b_1^{16}+b_1^{18},
\]
and, as a polynomial in $a$,
\[
\text{gcd}(\text{Res}(h_1',h_2;b_1),\text{Res}(T_2,T_3;k))=1+a,
\]
which is a contradiction.

When $T_{19}=0$, we have
\[
k=b_1^8+b_1^9+b_1^{10}+b_1^{12}+b_1^{16}+b_1^{17}+b_1^{18}+b_1^{19},
\]
and, as a polynomial in $a$,
\[
\text{gcd}(\text{Res}(h_1',h_2;b_1),\text{Res}(T_2,T_3;k))=1,
\]
which is a contradiction.

\medskip
{\bf Summary:} Only Case 2 is possible, and in that case we have $a\ne 1$ and
\begin{equation}\label{case2}
k=(1+a)^5(1+a^2+a^3)/a^2,\quad b_1=a/(1+a)^3,\quad E_3=0.
\end{equation}

\medskip
$3^\circ$ Since $E_3=0$, \eqref{2.10} -- \eqref{2.16} become
\begin{align}
\label{2.41}
&D_3^2=F_6,\\
\label{2.42}
&D_3E_2=F_5,\\
\label{2.43}
&D_2^2+D_3E_1+D_2E_2=F_4,\\
\label{2.44}
&D_3E_0+D_2E_1+D_1E_2=F_3,\\
\label{2.45}
&D_1^2+D_2E_0+D_1E_1+D_0E_2=F_2,\\
\label{2.46}
&D_1E_0+D_0E_1=F_1,\\
\label{2.47}
&D_0^2+D_0E_0=F_0.
\end{align}
We assume that $E_0\ne 0$ temporarily. By \eqref{2.44} -- \eqref{2.46}, 
\begin{align*}
&D_1=(D_0E_1+F_1)/E_0,\cr
&D_2=(D_1^2+D_1E_1+D_0E_2+F_2)/E_0,\cr
&D_3=(D_2E_1+D_1E_2+F_3)/E_0.
\end{align*}
By these substitutions, the remaining equations of \eqref{2.41} -- \eqref{2.47} become
\begin{align}
\label{2.48} 
&D_0^4 E_1^6+D_0^2 E_0^2 E_1^6+E_0^8
F_6+E_0^6 F_3^2+E_0^4 E_1^2
F_2^2\\
&+E_0^4 E_2^2 F_1^2+E_0^2
E_1^4 F_1^2+E_1^2 F_1^4=0, \cr
\label{2.49}
&D_0^2 E_1^3 E_2+D_0 E_0 E_1^3
E_2+E_0^4 F_5+E_0^3 E_2 F_3+E_0^2
E_1 E_2 F_2\\
&+E_0^2 E_2^2 F_1+E_0
E_1^2 E_2 F_1+E_1 E_2 F_1^2=0, \cr
\label{2.50}
&D_0^4 E_1^4+D_0^2 E_0^4 E_2^2+D_0^2
E_0^3 E_1^2 E_2+D_0 E_0^5
E_2^2+D_0 E_0^4 E_1^2 E_2\\
&+D_0E_0^3 E_1^4+E_0^6 F_4+E_0^5
E_1 F_3+E_0^5 E_2 F_2+E_0^4
E_1^2 F_2\cr
&+E_0^4 F_2^2+E_0^3
E_1^3 F_1+E_0^3 E_2 F_1^2+F_1^4=0, \cr
\label{2.51}
&D_0^2+D_0 E_0+F_0=0.
\end{align}
In fact, \eqref{2.48} -- \eqref{2.51} are linear combinations of \eqref{2.41} -- \eqref{2.47} with coefficients in $\f_2[D_0,\dots,D_3, E_0, E_1, E_2, F_0,\dots,F_6]$. Hence \eqref{2.48} -- \eqref{2.51} hold even if $E_0=0$. Thus we remove the assumption that $E_0\ne 0$. Subtracting \eqref{2.48} -- \eqref{2.50} by suitable multiples of \eqref{2.51} allows us to eliminate $D_0$, resulting in three equations in $E_i$'s and $F_j$'s:
\begin{align*}
L_1:=\,&E_0^8 F_6+E_0^6 F_3^2+E_0^4 E_1^2
F_2^2+E_0^4 E_2^2 F_1^2+E_0^2 E_1^4
F_1^2+E_1^6 F_0^2+E_1^2 F_1^4 =0,\cr
L_2:=\,&E_0^4 F_5+E_0^3 E_2 F_3+E_0^2
E_1 E_2 F_2+E_0^2 E_2^2
F_1+E_0 E_1^2 E_2 F_1+E_1^3
E_2 F_0\cr
&+E_1 E_2 F_1^2 =0,\cr
L_3:=\,&E_0^6 F_4+E_0^5 E_1 F_3+E_0^5
E_2 F_2+E_0^4 E_1^2 F_2+E_0^4
E_2^2 F_0+E_0^4 F_2^2+E_0^3
E_1^3 F_1\cr
&+E_0^3 E_1^2 E_2
F_0+E_0^3 E_2 F_1^2+E_0^2
E_1^4 F_0+E_1^4 F_0^2+F_1^4 =0.
\end{align*}
$L_1, L_2, L_3$ are polynomials in $a,b_1,k$. Using \eqref{case2}, we may express them in terms of $a$:
\begin{align*}
L_1\,&=\frac{a^{24}(1+a^2+a^3)^6(1+a+a^9)^2}{(1+a)^{40}},\cr
L_2\,&=\frac{a^{15}(1+a^2+a^3)^3(1+a+a^9)}{(1+a)^{25}},\cr
L_3\,&=\frac{a^{16}(1+a^2+a^3)^5(1+a^3+a^4)(1+a^4+a^5+a^6+a^7+a^{16}+a^{19}+a^{20}+a^{21})}{(1+a)^{36}}.
\end{align*}
It follows that $1+a^2+a^3=0$. But then, by \eqref{case2},
\[
k=\frac{(1+a)^5(1+a^2+a^3)}{a^2}=0,
\]
which is a contradiction.

This completes the proof of Claim~\ref{C3.1}.

\section{Completion of the Proof}

Form now on we assume that $b\in\f_q^*$. 

\medskip
$1^\circ$ We claim that if $f$ is a PP of $\f_{q^2}$, then $a=b$. Assume to the contrary that $a\ne b$. Choose $k\in\f_q$ such that $\text{Tr}_{q/2}(k)=1$ and let $z\in\f_{q^2}$ be such that $z^2+z+k=0$. We go through the computations in Section~3 again. However, since $b\in\f_q^*$, the computations are simpler. For \eqref{2.4.0} we have
\[
g\Bigl(\frac{X+z+1}{X+z}\Bigr)=\frac{A(X)}{B(X)},
\]
where
\begin{align*}
A(X)=\,& 1+az+az^2+az^3+bz^3+z^4+az^4+bz^4+ (a+az^2+bz^2 )X\cr
& +(a+az+bz )X^2 +(a+b )X^3+ (1+a+b )X^4,\cr
B(X)=\,& bz+bz^2+az^3+bz^3+z^4+az^4+bz^4 + (b+az^2+bz^2 )X\cr
& +(b+az+bz )X^2 +(a+b )X^3+ (1+a+b )X^4.
\end{align*}
Equation \eqref{2.4} becomes
\begin{equation}\label{4.0}
\frac{A(x)}{B(x)}=\frac{y+z+1}{y+z},
\end{equation}
which can be written as
\[
C_4x^4+C_3x^3+C_2(y)x^2+C_1(y)x+C_0(y)=0,
\]
where
\begin{align*}
C_4\,&=1+a+b\ (\ne 0),\cr
C_3\,&=a+b\ (\ne 0),\cr
C_2(Y)\,&=b+(a+b)Y,\cr
C_1(Y)\,&=b+ak+bk+(a+b)Y,\cr
C_0(Y)\,&=k+bk+k^2+ak^2+bk^2+(1+ak+bk)Y.
\end{align*}
We have
\[
C_1^2C_4+C_1C_2C_3+C_0C_3^2=b^2+a^2k+b^2k+(a+b)^2Y+(a+b)^2Y^2\ne0,
\]
hence $C_3(C_1^2C_4+C_1C_2C_3+C_0C_3^2)\ne0$. For \eqref{E0} -- \eqref{F6}, we have
\begin{align*}
E_0=\,& b^2+(a^2+b^2)k, \cr
E_1=\,& a^2+b^2,\cr
E_2=\,& a^2+b^2,\cr
E_3=\,& 0,\cr
F_0=\,& a^3b^3+b^4+(a^4b^2+a^2b^4)k+(a^4+a^5b+b^4+ab^5)k^2\cr
&+(a^6+a^4b^2+a^2b^4+b^6)k^3,\cr
F_1=\,& a^4b^2+a^2b^4+(a^6+a^4b^2+a^2b^4+b^6)k^2,\cr
F_2=\,& a^4+a^5b+a^4b^2+b^4+a^2b^4+ab^5+(a^6+a^4b^2+a^2b^4+b^6)k\cr
& +(a^6+a^4b^2+a^2b^4+b^6)k^2,\cr
F_3=\,& a^6+a^4b^2+a^2b^4+b^6,\cr
F_4=\,& a^4+a^6+a^5b+a^4b^2+b^4+a^2b^4+ab^5+b^6+(a^6+a^4b^2+a^2b^4+b^6)k,\cr
F_5=\,&a^6+a^4b^2+a^2b^4+b^6, \cr
F_6=\,& a^6+a^4b^2+a^2b^4+b^6.
\end{align*}
Since $E_3=0$, by \eqref{2.41} and \eqref{2.42}, we have $D_3^2=F_6$ and $D_3E_2=F_5$, which give $D_3=(a+b)^3$ and $(a+b)^5=(a+b)^6$. Then $(a+b)(a+b+1)=0$, which is a contradiction. Thus the claim is proved.

\medskip
$2^\circ$ Since $a=b$, we have $C_4=1$, $C_3=0$, $C_2=a$, $C_1=a$, $C_0=k+ak+k^2+Y$. Now \eqref{2.6} becomes
\begin{equation}\label{4.1} 
x^4+a^2x^2+a^2x+(k+ak+k^2+y)^2=0.
\end{equation}
We observe that
\begin{align*}
&\text{\eqref{4.1} has a unique solution $x$ in $\f_q$}\cr
\Leftrightarrow\ &\text{$X^3+a^2X+a^2$ is irreducible over $\f_q$ (Theorem~\ref{LW})}\cr
\Leftrightarrow\ &\text{$(X/a)^3+X/a+1/a$ has no root in $\f_q$}\cr
\Leftrightarrow\ &\text{$X^3+X+1/a$ has no root in $\f_q$}.
\end{align*}
This completes the proof for both necessity and sufficiency parts of the Theorem~\ref{T1.1}.



\section*{Appendix}

Gathered here are the lengthy intermediate computational results that were produced and used in the proof of Theorem~\ref{T1.1}.

In \eqref{2.24} -- \eqref{2.26},
\begin{align*}\tag{A1}\label{A1}
&h_1= a^{40}+a^{36}+b_1^3 a^{33}+b_1^4 a^{32}+b_1^3
a^{32}+b_1^4 a^{30}+b_1^3 a^{30}+b_1^4 a^{29}+b_1^6
a^{28}+b_1^4 a^{28}\\
&+b_1^3 a^{28}+b_1^6 a^{27}+b_1^5a^{27}+b_1^5 a^{26}+b_1^3 a^{26}+b_1^7 a^{25}+b_1^5a^{25}+b_1^4 a^{25}+b_1^8 a^{24}+b_1^7 a^{24}\cr
&+b_1^6a^{24}+b_1^4 a^{24}+b_1^3 a^{24}+a^{24}+b_1^5
a^{23}+b_1^8 a^{22}+b_1^7 a^{22}+b_1^6 a^{22}+b_1^5a^{22}+b_1^3 a^{22}\cr
&+b_1^7 a^{21}+b_1^6 a^{21}+b_1^5a^{21}+b_1^4 a^{21}+b_1^8 a^{20}+b_1^4 a^{20}+b_1^3a^{20}+a^{20}+b_1^9 a^{19}+b_1^8 a^{19}\cr
&+b_1^6a^{19}+b_1^5 a^{19}+b_1^9 a^{18}+b_1^5 a^{18}+b_1^3
a^{18}+b_1^9 a^{17}+b_1^7 a^{17}+b_1^5 a^{17}+b_1^4a^{17}+b_1^3 a^{17}\cr
&+b_1^8 a^{16}+b_1^7a^{16}+b_1^6 a^{16}+b_1^{10} a^{15}+b_1^5
a^{15}+b_1^{10} a^{14}+b_1^8 a^{14}+b_1^7a^{14}+b_1^6 a^{14}+b_1^5a^{14}\cr
&+b_1^4a^{14}+b_1^{10} a^{13}+b_1^7 a^{13}+b_1^6
a^{13}+b_1^5 a^{13}+b_1^{10} a^{12}+b_1^8a^{12}+b_1^6 a^{12}+b_1^{10} a^{11}+b_1^9a^{11}\cr
&+b_1^8 a^{11}+b_1^{10} a^{10}+b_1^9a^{10}+b_1^{10} a^9+b_1^9 a^9+b_1^{10} a^8+b_1^{40}k^{20}+b_1^{36} k^{18}+b_1^{36} k^{17}+(b_1^{36}\cr
&+ab_1^{35}+b_1^{35}+a^8 b_1^{32}+a^4
b_1^{32}) k^{16}+(b_1^{34}+b_1^{33})k^{15}+(b_1^{34}+a^2 b_1^{32}+ab_1^{32}+b_1^{32}\cr
&+a^2 b_1^{31}+b_1^{31})k^{14}+(a b_1^{32}+a b_1^{31}+a^4 b_1^{30}+a^4b_1^{29}+b_1^{29}) k^{13}+(b_1^{32}+ab_1^{31}+b_1^{31}\cr
&+a^4 b_1^{30}+a^3b_1^{30}+a^3 b_1^{29}+a^2 b_1^{29}+ab_1^{29}+b_1^{29}+a^6 b_1^{28}+a^5 b_1^{28}+a^4b_1^{28}+a b_1^{28}+b_1^{28}\cr
&+a^6 b_1^{27}+a^4b_1^{27}+a^2 b_1^{27}+b_1^{27}+b_1^{24})k^{12}+(b_1^{29}+a^2 b_1^{28}+a b_1^{28}+a^2b_1^{27}+a b_1^{27}+b_1^{27}\cr
&+a^8 b_1^{26}+a^8b_1^{25}+b_1^{25}) k^{11}+(a^2b_1^{28}+b_1^{28}+a^2 b_1^{27}+b_1^{27}+a^8b_1^{26}+a^4 b_1^{26}+a^2 b_1^{26}\cr
&+ab_1^{26}+b_1^{26}+a^4 b_1^{25}+a^3 b_1^{25}+a^2
b_1^{25}+a b_1^{25}+a^{10} b_1^{24}+a^9 b_1^{24}+a^8b_1^{24}+a b_1^{24}+b_1^{24}\cr
&+a^{10} b_1^{23}+a^8b_1^{23}+a^2 b_1^{23}+b_1^{23}+b_1^{20})
k^{10}+(a^2 b_1^{26}+a b_1^{26}+b_1^{26}+a^4b_1^{25}+a^3 b_1^{25}\cr
&+ab_1^{25}+b_1^{25}+a^9b_1^{24}+a^6 b_1^{24}+a^5 b_1^{24}+a^2
b_1^{24}+a b_1^{24}+a^9 b_1^{23}+a^6 b_1^{23}+a^5
b_1^{23}+a^4 b_1^{23}\cr
&+a b_1^{23}+a^{12}b_1^{22}+a^{12} b_1^{21}+a^8 b_1^{21}+a^4
b_1^{21}+b_1^{21}+b_1^{20}) k^9+(a^3b_1^{25}+a^2 b_1^{25}+a b_1^{25}\cr
&+b_1^{25}+a^8b_1^{24}+a^6 b_1^{24}+a^4 b_1^{24}+a^3 b_1^{24}+a^2
b_1^{24}+a b_1^{24}+b_1^{24}+a^9 b_1^{23}+a^8b_1^{23}+a^6 b_1^{23}\cr
&+a^5 b_1^{23}+a^2 b_1^{23}+ab_1^{23}+a^{12} b_1^{22}+a^{11} b_1^{22}+a^8
b_1^{22}+a^6 b_1^{22}+a^5 b_1^{22}+a^3b_1^{22}+a^{11} b_1^{21}\cr
&+a^{10} b_1^{21}+a^9b_1^{21}+a^7 b_1^{21}+a^6 b_1^{21}+a^5 b_1^{21}+a^3
b_1^{21}+a^2 b_1^{21}+a b_1^{21}+b_1^{21}+a^{14}b_1^{20}\cr
&+a^{13} b_1^{20}+a^{12} b_1^{20}+a^9b_1^{20}+a^8 b_1^{20}+a^5 b_1^{20}+a^4 b_1^{20}+ab_1^{20}+a^{14} b_1^{19}+a^{12} b_1^{19}+a^{10}b_1^{19}\cr
&+a^8 b_1^{19}+a^6 b_1^{19}+a^4 b_1^{19}+a^2b_1^{19}+a b_1^{19}+a^8 b_1^{16}+a^4 b_1^{16})k^8+(a^2 b_1^{23}+a b_1^{23}+b_1^{23}\cr
&+a^4b_1^{22}+a^4 b_1^{21}+a^2 b_1^{20}+a
b_1^{20}+a^2 b_1^{19}+a b_1^{19}+b_1^{19}+a^{16}
b_1^{18}+b_1^{18}+a^{16} b_1^{17})k^7\cr
&+(b_1^{22}+a^4 b_1^{21}+a^3 b_1^{21}+a^2
b_1^{21}+a b_1^{21}+a^6 b_1^{20}+a^5 b_1^{20}+a^4
b_1^{20}+a^2 b_1^{20}+a b_1^{20}+b_1^{20}\cr
&+a^6b_1^{19}+a^4 b_1^{19}+a^{16} b_1^{18}+a^4
b_1^{18}+a^2 b_1^{18}+a b_1^{18}+a^4 b_1^{17}+a^3
b_1^{17}+a^2 b_1^{17}+a b_1^{17}\cr
&+a^{18}b_1^{16}+a^{17} b_1^{16}+a^{16} b_1^{16}+a^2b_1^{16}+a^{18} b_1^{15}+a^{16} b_1^{15})k^6+(b_1^{14} a^{20}+b_1^{13} a^{20}+b_1^{16}
a^{17}\cr
&+b_1^{15} a^{17}+b_1^{13} a^{16}+b_1^{18}a^8+b_1^{17} a^8+b_1^{19} a^6+b_1^{16} a^6+b_1^{15}a^6+b_1^{20} a^5+b_1^{16} a^5+b_1^{15} a^5\cr
&+b_1^{20}a^4+b_1^{19} a^4+b_1^{17} a^4+b_1^{15} a^4+b_1^{14}
a^4+b_1^{20} a^3+b_1^{17} a^3+b_1^{18}a^2+b_1^{16} a^2+b_1^{19} a\cr
&+b_1^{18} a+b_1^{17}a+b_1^{20}+b_1^{18}) k^5+(b_1^8a^{32}+b_1^{12} a^{22}+b_1^{11} a^{22}+b_1^{12}a^{21}+b_1^{14} a^{20}+b_1^{12} a^{20}\cr
&+b_1^{11}a^{20}+b_1^{14} a^{19}+b_1^{13} a^{19}+b_1^{13}
a^{18}+b_1^{11} a^{18}+b_1^{15} a^{17}+b_1^{13}
a^{17}+b_1^{12} a^{17}+b_1^{16} a^{16}\cr
&+b_1^{15}a^{16}+b_1^{13} a^{16}+b_1^{12} a^{16}+b_1^{11}
a^{16}+b_1^8 a^{16}+b_1^{16} a^{10}+b_1^{15}
a^{10}+b_1^{16} a^9+b_1^{17} a^8\cr
&+b_1^{15}a^8+b_1^{14} a^8+b_1^{13} a^8+b_1^{18} a^7+b_1^{13}
a^7+b_1^{18} a^6+b_1^{16} a^6+b_1^{15} a^6+b_1^{14}
a^6+b_1^{13} a^6\cr
&+b_1^{12} a^6+b_1^{18} a^5+b_1^{15}a^5+b_1^{14} a^5+b_1^{13} a^5+b_1^{17} a^4+b_1^{14}a^4+b_1^{18} a^3+b_1^{16} a^3+b_1^{18} a^2\cr
&+b_1^{16}a^2+b_1^{18} a) k^4+(b_1^{10} a^{24}+b_1^9
a^{24}+b_1^{12} a^{18}+b_1^{11} a^{18}+b_1^{12}
a^{17}+b_1^{11} a^{17}+b_1^{13} a^{16}\cr
&+b_1^{11}a^{16}+b_1^9 a^{16}+b_1^{14} a^{12}+b_1^{13}
a^{12}+b_1^{15} a^{10}+b_1^{12} a^{10}+b_1^{11}
a^{10}+b_1^{15} a^9+b_1^{12} a^9\cr
&+b_1^{11}a^9+b_1^{15} a^8+b_1^{11} a^8+b_1^{10} a^8+b_1^{16}
a^6+b_1^{15} a^6+b_1^{16} a^5+b_1^{15} a^5+b_1^{16}
a^4+b_1^{15} a^4\cr
&+b_1^{14} a^4+b_1^{17} a^3+b_1^{16}a^3+b_1^{17} a) k^3+(b_1^4 a^{32}+b_1^8
a^{26}+b_1^7 a^{26}+b_1^8 a^{25}+b_1^{10}a^{24}\cr
&+b_1^8 a^{24}+b_1^7 a^{24}+b_1^{10}a^{20}+b_1^9 a^{20}+b_1^9 a^{19}+b_1^{12}
a^{18}+b_1^{11} a^{18}+b_1^{10} a^{18}+b_1^9
a^{18}+b_1^7 a^{18}\cr
&+b_1^{10} a^{17}+b_1^9a^{17}+b_1^8 a^{17}+b_1^{12} a^{16}+b_1^{11}
a^{16}+b_1^{10} a^{16}+b_1^8 a^{16}+b_1^7a^{16}+b_1^4 a^{16}+b_1^{12} a^{14}\cr
&+b_1^{11}a^{14}+b_1^{12} a^{13}+b_1^{13} a^{12}+b_1^{12}
a^{12}+b_1^{11} a^{12}+b_1^{10} a^{12}+b_1^9a^{12}+b_1^{13} a^{11}+b_1^9 a^{11}\cr
&+b_1^{13}a^{10}+b_1^{12} a^{10}+b_1^{10} a^{10}+b_1^9
a^{10}+b_1^8 a^{10}+b_1^{13} a^9+b_1^{12}
a^9+b_1^{10} a^9+b_1^9 a^9+b_1^{13} a^8\cr
&+b_1^{13}a^7+b_1^{13} a^6+b_1^{12} a^6+b_1^{13} a^5)
k^2+(b_1^4 a^{32}+b_1^6 a^{28}+b_1^5
a^{28}+b_1^8 a^{25}+b_1^7 a^{25}+b_1^5 a^{24}\cr
&+b_1^8a^{22}+b_1^7 a^{22}+b_1^8 a^{21}+b_1^7 a^{21}+b_1^9
a^{20}+b_1^7 a^{20}+b_1^5 a^{20}+b_1^9
a^{19}+b_1^{10} a^{18}+b_1^8 a^{18}\cr
&+b_1^{10}a^{17}+b_1^9 a^{17}+b_1^8 a^{17}+b_1^7 a^{17}+b_1^5
a^{16}+b_1^4 a^{16}+b_1^{11} a^{14}+b_1^8
a^{14}+b_1^7 a^{14}+b_1^{12} a^{13}\cr
&+b_1^8a^{13}+b_1^7 a^{13}+b_1^{12} a^{12}+b_1^{11}
a^{12}+b_1^9 a^{12}+b_1^7 a^{12}+b_1^6
a^{12}+b_1^{12} a^{11}+b_1^9 a^{11}+b_1^{12}a^{10}\cr
&+b_1^{11} a^{10}+b_1^{10} a^{10}+b_1^8
a^{10}+b_1^{12} a^9+b_1^{10} a^9+b_1^9 a^9+b_1^{12}
a^8+b_1^{11} a^8+b_1^{12} a^7+b_1^{12} a^6) k, 
\end{align*}
\begin{align*}\tag{A2}\label{A2}
&h_2= a^{20}+a^{18}+a^{17} b_1+a^{16} b_1^4 k^2+a^{16} b_1^2
k+a^{15} b_1^3+a^{13} b_1^5 k+a^{12} b_1^4+a^{12}+a^{11}b_1^7 k^2\\
&+a^{10} b_1^4+a^{10}+a^9 b_1^9 k^3+a^9b_1^5+a^9 b_1+a^8 b_1^8 k^2+a^8 b_1^6 k+a^8b_1^4 k^2+a^8 b_1^4+a^8 b_1^2 k\cr
&+a^7 b_1^{11}k^4+a^7 b_1^7 k+a^7 b_1^3+a^6 b_1^4+a^5 b_1^{13}
k^5+a^5 b_1^9 k^2+a^5 b_1^5 k+a^5 b_1^5+a^4b_1^{16} k^8\cr
&+a^4 b_1^{12} k^4+a^4 b_1^8 k^4+a^4b_1^8 k^2+a^4 b_1^6 k+a^3 b_1^{15} k^6+a^3 b_1^{11}k^3+a^3 b_1^7 k^2+a^3 b_1^7 k\cr
&+a^2 b_1^{16} k^8+a^2b_1^{12} k^4+a^2 b_1^8 k^4+a b_1^{17} k^8+a
b_1^{17} k^7+a b_1^9 k^4+a b_1^9 k^3+a b_1^9
k^2+b_1^{20} k^{10}\cr
&+b_1^{18} k^9+b_1^{16}k^6+b_1^{14} k^5+b_1^{12} k^6+b_1^{10} k^5,
\end{align*} 
\begin{align*}\tag{A3}\label{A3}
&h_3= a^{50}+b_1^2 a^{48}+a^{48}+b_1^2 a^{44}+b_1^4
a^{42}+b_1^2 a^{42}+a^{42}+b_1^6 a^{40}+a^{40}+b_1^8a^{38}+b_1^6 a^{38}\\
&+b_1^4 a^{38}+b_1^2a^{38}+b_1^6 a^{36}+b_1^6 a^{34}+a^{34}+b_1^{10}
a^{32}+b_1^8 a^{32}+b_1^2 a^{32}+a^{32}+b_1^{10}a^{28}\cr
&+b_1^2 a^{28}+b_1^{12} a^{26}+b_1^{10}a^{26}+b_1^4 a^{26}+b_1^2 a^{26}+a^{26}+b_1^{10}a^{24}+b_1^6 a^{24}+a^{24}+b_1^{12} a^{22}\cr
&+b_1^{10}a^{22}+b_1^8 a^{22}+b_1^6 a^{22}+b_1^4
a^{22}+b_1^2 a^{22}+b_1^{10} a^{20}+b_1^6a^{20}+b_1^{12} a^{18}+b_1^{10} a^{18}+b_1^6a^{18}\cr
&+b_1^8 a^{16}+b_1^{12} a^{14}+b_1^{10}a^{14}+(b_1^{50}+a^2 b_1^{48}+b_1^{48})
k^{24}+(a b_1^{45}+b_1^{45}+a^2 b_1^{44})k^{21}+(a^2b_1^{44}\cr
&+b_1^{44}+b_1^{43}+a^4b_1^{42}+a^2 b_1^{42}+b_1^{42}+a^2
b_1^{40}+b_1^{40}) k^{20}+(a b_1^{43}+a^2b_1^{42}+a b_1^{42}+b_1^{42}\cr
&+a^2 b_1^{41}+ab_1^{41}+a^2 b_1^{40}) k^{19}+(a^2 b_1^{44}+a
b_1^{43}+a^4 b_1^{42}+a b_1^{42}+b_1^{42}+a^4b_1^{40}+a^2 b_1^{40}\cr
&+a^4 b_1^{39}+a^4 b_1^{38}+a^2b_1^{38}) k^{18}+(a b_1^{41}+b_1^{41}+a^4
b_1^{40}+a^3 b_1^{40}+a^2 b_1^{40}+a^5 b_1^{39}+a^3b_1^{39}\cr
&+a^6 b_1^{38}+a^5 b_1^{38}+a^9 b_1^{37}+a^8b_1^{37}+a^6 b_1^{37}+a^5 b_1^{37}+a^{10}b_1^{36}+a^6 b_1^{36}) k^{17}+(b_1^{42}+a^6b_1^{40}\cr
&+b_1^{40}+b_1^{39}+a^8 b_1^{38}+a^6b_1^{38}+a^4 b_1^{38}+a^2 b_1^{38}+b_1^{38}+a^{10}b_1^{36}+a^6 b_1^{36}+a^{16} b_1^{34}+a^{12}b_1^{34}\cr
&+a^{10} b_1^{34}+a^6 b_1^{34}+b_1^{34}+a^{18}b_1^{32}+a^{16} b_1^{32}+a^{10} b_1^{32}+a^8b_1^{32}+a^2 b_1^{32}+b_1^{32}) k^{16}+(ab_1^{39}\cr
&+a^2 b_1^{38}+a^2 b_1^{37}+a b_1^{37}+a^2b_1^{36}) k^{15}+(a^3 b_1^{37}+a^3 b_1^{36}+a^2b_1^{36}+b_1^{36}+a^4 b_1^{35}+a^4 b_1^{34})k^{14}\cr
&+(b_1^{36}+a^6 b_1^{34}+a^5 b_1^{34}+ab_1^{34}+b_1^{34}+a^6 b_1^{33}+a^5 b_1^{33}+a^6b_1^{32}+a b_1^{29}+b_1^{29}+a^2 b_1^{28})k^{13}\cr
&+(a^4 b_1^{34}+a^2 b_1^{34}+a^3 b_1^{33}+a^3b_1^{32}+a^2 b_1^{32}+a^8 b_1^{31}+a^8 b_1^{30}+a^2b_1^{28}+b_1^{28}+b_1^{27}+a^4 b_1^{26}\cr
&+a^2b_1^{26}+b_1^{26}+a^2 b_1^{24}+b_1^{24})k^{12}+(a^3 b_1^{33}+a^9 b_1^{31}+a^5 b_1^{31}+a^3b_1^{31}+a^{10} b_1^{30}+a^5 b_1^{30}\cr
&+a^4b_1^{30}+a b_1^{30}+b_1^{30}+a^{10} b_1^{29}+a^9b_1^{29}+a^{10} b_1^{28}+a b_1^{27}+a^2 b_1^{26}+ab_1^{26}+b_1^{26}+a^2 b_1^{25}\cr
&+a b_1^{25}+a^2b_1^{24}) k^{11}+(a^6 b_1^{32}+a^5b_1^{31}+a b_1^{31}+a^5 b_1^{30}+a^2 b_1^{30}+ab_1^{30}+a^{11} b_1^{29}+a^7 b_1^{29}\cr
&+a^3b_1^{29}+a^{11} b_1^{28}+a^{10} b_1^{28}+a^8b_1^{28}+a^7 b_1^{28}+a^4 b_1^{28}+a^3 b_1^{28}+a^2b_1^{28}+b_1^{28}+a^{12} b_1^{27}\cr
&+a b_1^{27}+a^{12}b_1^{26}+a^4 b_1^{26}+a b_1^{26}+b_1^{26}+a^4b_1^{24}+a^2 b_1^{24}+a^4 b_1^{23}+a^4 b_1^{22}+a^2b_1^{22}) k^{10}\cr
&+(a^5 b_1^{29}+a^4b_1^{29}+a b_1^{29}+b_1^{29}+a^8 b_1^{28}+a^7
b_1^{28}+a^6 b_1^{28}+a^4 b_1^{28}+a^3 b_1^{28}+a^2b_1^{28}+b_1^{28}\cr
&+a^9 b_1^{27}+a^7 b_1^{27}+a^5b_1^{27}+a^3 b_1^{27}+a^{14} b_1^{26}+a^{13}
b_1^{26}+a b_1^{26}+b_1^{26}+a^{14} b_1^{25}+a^{13}b_1^{25}\cr
&+a b_1^{25}+b_1^{25}+a^{14} b_1^{24}+a^4b_1^{24}+a^3 b_1^{24}+a^2 b_1^{24}+a^5 b_1^{23}+a^3b_1^{23}+a^6 b_1^{22}+a^5 b_1^{22}\cr
&+a^9 b_1^{21}+a^8b_1^{21}+a^6 b_1^{21}+a^5 b_1^{21}+a^{10}
b_1^{20}+a^6 b_1^{20}) k^9+(b_1^{16}a^{34}+b_1^{18} a^{32}+b_1^{16} a^{32}\cr
&+b_1^{23}a^{16}+b_1^{22} a^{16}+b_1^{26} a^{12}+b_1^{18}a^{12}+b_1^{28} a^{10}+b_1^{26} a^{10}+b_1^{20}a^{10}+b_1^{18} a^{10}+b_1^{16} a^{10}\cr
&+b_1^{26}a^8+b_1^{22} a^8+b_1^{16} a^8+b_1^{28} a^6+b_1^{26}a^6+b_1^{24} a^6+b_1^{22} a^6+b_1^{20} a^6+b_1^{18}a^6+b_1^{27} a^4\cr
&+b_1^{26} a^4+b_1^{22} a^4+b_1^{25}a^3+b_1^{24} a^3+b_1^{24} a^2+b_1^{22}
a^2+b_1^{27}+b_1^{26}+b_1^{24}+b_1^{23}+b_1^{22}) k^8\cr
&+(a^5 b_1^{27}+a b_1^{27}+a^6b_1^{26}+a^2 b_1^{26}+a^6 b_1^{25}+a^5 b_1^{25}+a^3b_1^{25}+a^2 b_1^{25}+a b_1^{25}+a^6 b_1^{24}\cr
&+a^2b_1^{24}+a^{17} b_1^{23}+a^5 b_1^{23}+a^3 b_1^{23}+ab_1^{23}+a^{18} b_1^{22}+a^5 b_1^{22}+a^4b_1^{22}+a^2 b_1^{22}+a b_1^{22}\cr
&+b_1^{22}+a^{18}b_1^{21}+a^{17} b_1^{21}+a^2 b_1^{21}+ab_1^{21}+a^{18} b_1^{20}+a^2 b_1^{20})k^7+(a^4 b_1^{26}+b_1^{26}+a^7 b_1^{25}\cr
&+a^3b_1^{25}+a^7 b_1^{24}+a^6 b_1^{24}+a^4 b_1^{24}+a^3b_1^{24}+a^2 b_1^{24}+b_1^{24}+a^8 b_1^{23}+a^5b_1^{23}+a^4 b_1^{23}+a b_1^{23}\cr
&+a^8 b_1^{22}+a^5b_1^{22}+a^4 b_1^{22}+a b_1^{22}+a^{19}b_1^{21}+a^7 b_1^{21}+a^{19} b_1^{20}+a^{18}b_1^{20}+a^{16} b_1^{20}+a^7 b_1^{20}\cr
&+a^4b_1^{20}+a^2 b_1^{20}+a^{20} b_1^{19}+a^4b_1^{19}+a^{20} b_1^{18}+a^4 b_1^{18})k^6+(b_1^{12} a^{34}+b_1^{13} a^{33}+b_1^{13}a^{32}\cr
&+b_1^{18} a^{22}+b_1^{17} a^{22}+b_1^{16}a^{22}+b_1^{18} a^{21}+b_1^{17} a^{21}+b_1^{12}a^{18}+b_1^{18} a^{17}+b_1^{13} a^{17}+b_1^{20}a^{16}\cr
&+b_1^{18} a^{16}+b_1^{13} a^{16}+b_1^{22}a^{10}+b_1^{21} a^{10}+b_1^{22} a^9+b_1^{19}a^9+b_1^{18} a^9+b_1^{21} a^8+b_1^{18} a^8+b_1^{20}a^7\cr
&+b_1^{19} a^7+b_1^{22} a^6+b_1^{21} a^6+b_1^{18}a^6+b_1^{17} a^6+b_1^{16} a^6+b_1^{19} a^5+b_1^{18}a^5+b_1^{17} a^5+b_1^{24} a^4\cr
&+b_1^{22} a^4+b_1^{21}a^4+b_1^{20} a^4+b_1^{20} a^3+b_1^{19} a^3+b_1^{22}
a+b_1^{24}+b_1^{22}) k^5+(b_1^{10}a^{36}+b_1^{12} a^{34}\cr
&+b_1^{10} a^{34}+b_1^8a^{34}+b_1^{12} a^{32}+b_1^{11} a^{32}+b_1^{10}
a^{32}+b_1^8 a^{32}+b_1^{15} a^{24}+b_1^{14}a^{24}+b_1^{18} a^{20}\cr
&+b_1^{10} a^{20}+b_1^{17}a^{19}+b_1^{16} a^{19}+b_1^{18} a^{18}+b_1^{16}
a^{18}+b_1^{12} a^{18}+b_1^{10} a^{18}+b_1^8a^{18}+b_1^{12} a^{16}\cr
&+b_1^{11} a^{16}+b_1^{10}a^{16}+b_1^8 a^{16}+b_1^{19} a^{12}+b_1^{17}
a^{11}+b_1^{16} a^{11}+b_1^{20} a^{10}+b_1^{18}a^{10}+b_1^{16} a^{10}\cr
&+b_1^{22} a^8+b_1^{20}a^8+b_1^{15} a^8+b_1^{14} a^8+b_1^{21} a^7+b_1^{20}
a^7+b_1^{22} a^6+b_1^{22} a^4+b_1^{20} a^4+b_1^{19}a^4\cr
&+b_1^{18} a^4+b_1^{21} a^3+b_1^{20} a^3+b_1^{22}a^2+b_1^{20} a^2) k^4+(b_1^{10} a^{34}+b_1^9a^{34}+b_1^8 a^{34}+b_1^{11} a^{33}\cr
&+b_1^{10}a^{33}+b_1^9 a^{33}+b_1^{10} a^{32}+b_1^{14}
a^{26}+b_1^{13} a^{26}+b_1^{12} a^{26}+b_1^{15}a^{25}+b_1^{13} a^{25}+b_1^{15} a^{21}\cr
&+b_1^{14}a^{21}+b_1^{14} a^{20}+b_1^{17} a^{19}+b_1^{15}a^{19}+b_1^{10} a^{18}+b_1^9 a^{18}+b_1^8a^{18}+b_1^{14} a^{17}+b_1^{11} a^{17}\cr
&+b_1^{10}a^{17}+b_1^9 a^{17}+b_1^{14} a^{16}+b_1^{10}a^{16}+b_1^{18} a^{14}+b_1^{17} a^{14}+b_1^{16}a^{14}+b_1^{19} a^{13}+b_1^{17} a^{13}\cr
&+b_1^{15}a^{13}+b_1^{14} a^{13}+b_1^{14} a^{12}+b_1^{17}a^{11}+b_1^{15} a^{11}+b_1^{14} a^{10}+b_1^{13}a^{10}+b_1^{12} a^{10}+b_1^{19} a^9\cr
&+b_1^{15}a^9+b_1^{14} a^9+b_1^{13} a^9+b_1^{14} a^8+b_1^{21}a^7+b_1^{19} a^7+b_1^{18} a^6+b_1^{17} a^6+b_1^{16}a^6+b_1^{17} a^5\cr
&+b_1^{21} a^3+b_1^{19} a^3)k^3+(b_1^{10} a^{36}+b_1^8 a^{36}+b_1^7
a^{36}+b_1^6 a^{36}+b_1^{12} a^{34}+b_1^8a^{34}+b_1^6 a^{34}\cr
&+b_1^{11} a^{33}+b_1^{10}a^{33}+b_1^{10} a^{32}+b_1^{11} a^{28}+b_1^{10}
a^{28}+b_1^{13} a^{27}+b_1^{12} a^{27}+b_1^{12}a^{26}+b_1^{12} a^{24}\cr
&+b_1^{13} a^{23}+b_1^{12}a^{23}+b_1^{16} a^{22}+b_1^{15} a^{21}+b_1^{14}
a^{21}+b_1^{12} a^{20}+b_1^{10} a^{20}+b_1^8a^{20}+b_1^7 a^{20}\cr
&+b_1^6 a^{20}+b_1^{13}a^{19}+b_1^{12} a^{19}+b_1^{14} a^{18}+b_1^{12}
a^{18}+b_1^8 a^{18}+b_1^6 a^{18}+b_1^{15}a^{17}+b_1^{14} a^{17}\cr
&+b_1^{11} a^{17}+b_1^{10}a^{17}+b_1^{15} a^{16}+b_1^{14} a^{16}+b_1^{12}
a^{16}+b_1^{10} a^{16}+b_1^{17} a^{15}+b_1^{16}a^{15}+b_1^{13} a^{15}\cr
&+b_1^{12} a^{15}+b_1^{16}a^{14}+b_1^{15} a^{13}+b_1^{14} a^{13}+b_1^{12}
a^{12}+b_1^{11} a^{12}+b_1^{10} a^{12}+b_1^{14}a^{10}+b_1^{12} a^{10}\cr
&+b_1^{15} a^9+b_1^{14}a^9+b_1^{15} a^8+b_1^{14} a^8+b_1^{17} a^7+b_1^{16}
a^7) k^2+(b_1^4 a^{42}+b_1^5 a^{41}+b_1^5a^{40}+b_1^6 a^{38}\cr
&+b_1^5 a^{38}+b_1^4 a^{38}+b_1^7a^{37}+b_1^6 a^{37}+b_1^5 a^{37}+b_1^8 a^{36}+b_1^8a^{35}+b_1^7 a^{35}+b_1^8 a^{34}+b_1^9 a^{33}\cr
&+b_1^9a^{32}+b_1^{10} a^{30}+b_1^9 a^{30}+b_1^8a^{30}+b_1^{10} a^{29}+b_1^9 a^{29}+b_1^4a^{26}+b_1^{11} a^{25}+b_1^5 a^{25}+b_1^{12}a^{24}\cr
&+b_1^5 a^{24}+b_1^{12} a^{23}+b_1^{11}a^{23}+b_1^{12} a^{22}+b_1^6 a^{22}+b_1^5
a^{22}+b_1^4 a^{22}+b_1^{13} a^{21}+b_1^{11}a^{21}+b_1^7 a^{21}\cr
&+b_1^6 a^{21}+b_1^5a^{21}+b_1^{13} a^{20}+b_1^{12} a^{20}+b_1^8a^{20}+b_1^{12} a^{19}+b_1^{11} a^{19}+b_1^8a^{19}+b_1^7 a^{19}+b_1^{14} a^{18}\cr
&+b_1^{13}a^{18}+b_1^{12} a^{18}+b_1^8 a^{18}+b_1^{14}a^{17}+b_1^{13} a^{17}+b_1^{11} a^{17}+b_1^9a^{17}+b_1^{12} a^{16}+b_1^9 a^{16}+b_1^{12}a^{15}\cr
&+b_1^{11} a^{15}+b_1^{12} a^{14}+b_1^{10}a^{14}+b_1^9 a^{14}+b_1^8 a^{14}+b_1^{13}a^{13}+b_1^{11} a^{13}+b_1^{10} a^{13}+b_1^9a^{13}+b_1^{13} a^{12}\cr
&+b_1^{12} a^{12}+b_1^{12}
a^{11}+b_1^{11} a^{11}+b_1^{14} a^{10}+b_1^{13}
a^{10}+b_1^{12} a^{10}+b_1^{14} a^9+b_1^{13} a^9) k.
\end{align*}
In \eqref{eq3.36},
\begin{align*}\tag{A4}\label{A4}
&h_1'= a^{22}+b_1^4 a^{20}+b_1^2 a^{20}+a^{20}+b_1 a^{19}+b_1^5a^{17}+b_1^3 a^{17}+b_1^4 a^{16}+b_1^6a^{14}+a^{14}+b_1^6 a^{12}\\
&+b_1^2 a^{12}+a^{12}+b_1a^{11}+b_1^4 a^{10}+b_1^7 a^9+b_1^3 a^9+b_1^6
a^8+b_1^6 a^6+b_1^4 a^6+b_1^7 a^5+b_1^5a^5\cr
&+(a^3 b_1^{19}+a^4 b_1^{18}+a^2 b_1^{18})k^9+(a^4 b_1^{20}+a^2 b_1^{20}+b_1^{20}+a^4b_1^{18}+a^3 b_1^{17}+a^6 b_1^{16}+a^4b_1^{16}) k^8\cr
&+(a^3 b_1^{19}+a^4 b_1^{18}+a^2b_1^{18}+b_1^{18}+a^3 b_1^{17}+a b_1^{17})
k^7+(b_1^{18}+a^5 b_1^{17}+a b_1^{17}+a^6b_1^{16}+a^4 b_1^{16}\cr
&+a^2 b_1^{16}+b_1^{16}+a^5b_1^{15}+a b_1^{15}) k^6+(b_1^{16}+a^7
b_1^{15}+a^8 b_1^{14}+b_1^{14}+a^7 b_1^{13}+a^5b_1^{13}+a^3 b_1^{13}\cr
&+a b_1^{13}+a^3 b_1^{11}+a^4b_1^{10}+a^2 b_1^{10}) k^5+(a^6 b_1^{14}+a^2
b_1^{14}+a^9 b_1^{13}+a b_1^{13}+a^{10}b_1^{12}+a^4 b_1^{12}\cr
&+a^9 b_1^{11}+a b_1^{11}+a^4b_1^{10}+a^3 b_1^9+a^6 b_1^8+a^4 b_1^8)
k^4+(a^3 b_1^{13}+a^6 b_1^{12}+a^4 b_1^{12}+a^2b_1^{12}\cr
&+a^{11} b_1^{11}+a^3 b_1^{11}+ab_1^{11}+a^{12} b_1^{10}+a^{10} b_1^{10}+a^8
b_1^{10}+a^6 b_1^{10}+a^{11} b_1^9+a^9 b_1^9)k^3+(b_1^8 a^{14}\cr
&+b_1^9 a^{13}+b_1^7a^{13}+b_1^8 a^{12}+b_1^8 a^{10}+b_1^9 a^9+b_1^7
a^9+b_1^{11} a^5+b_1^9 a^5+b_1^8 a^4+b_1^{11}a+b_1^9 a) k^2\cr
&+(b_1^2 a^{20}+b_1^3a^{19}+b_1^2 a^{18}+b_1^6 a^{16}+b_1^7 a^{15}+b_1^5
a^{15}+b_1^5 a^{13}+b_1^2 a^{12}+b_1^5 a^{11}+b_1^3a^{11}\cr
&+b_1^8 a^{10}+b_1^6 a^{10}+b_1^2 a^{10}+b_1^5
a^9+b_1^8 a^8+b_1^9 a^7+b_1^8 a^6+b_1^6
a^6+b_1^8 a^4+b_1^6 a^4+b_1^9 a^3\cr
&+b_1^7 a^3) k, 
\end{align*}
\begin{align*}\tag{A5}\label{A5}
&h_3'= b_1^2 a^{24}+a^{24}+b_1^3 a^{23}+b_1^4 a^{22}+b_1^2a^{22}+b_1^5 a^{21}+b_1^4 a^{20}+b_1^2a^{20}+a^{20}+b_1^5 a^{19}\\
&+b_1^3 a^{19}+b_1^4a^{18}+b_1^7 a^{17}+b_1^5 a^{17}+b_1^8 a^{16}+b_1^6
a^{16}+b_1^2 a^{16}+a^{16}+b_1^5 a^{15}+b_1^3a^{15}+b_1^8 a^{14}\cr
&+b_1^6 a^{14}+b_1^4 a^{14}+b_1^2a^{14}+b_1^9 a^{13}+b_1^5 a^{13}+b_1^8 a^{12}+b_1^4a^{12}+b_1^2 a^{12}+a^{12}+b_1^7 a^{11}+b_1^5a^{11}\cr
&+b_1^3 a^{11}+b_1^{10} a^{10}+b_1^6a^{10}+b_1^4 a^{10}+b_1^5 a^9+b_1^6 a^8+b_1^7a^7+b_1^5 a^7+b_1^{10} a^6+b_1^8 a^6+b_1^9a^5\cr
&+b_1^7 a^5+(a^4 b_1^{24}+a^3 b_1^{23}+a^4b_1^{22}+a^5 b_1^{21}+a^3 b_1^{21}+a^4
b_1^{20}+a^5 b_1^{19}+a b_1^{19}+a^6 b_1^{18}\cr
&+a^2b_1^{18}) k^9+(a^3 b_1^{23}+a^6 b_1^{22}+a^2b_1^{22}+b_1^{22}+a^5 b_1^{21}+a^3 b_1^{21}+ab_1^{21}+a^2 b_1^{20}+a^5 b_1^{19}\cr
&+a^3 b_1^{19}+a^8b_1^{18}+a^6 b_1^{18}+a^2 b_1^{18}+a^8 b_1^{16}+a^4
b_1^{16}) k^8+(a^2 b_1^{22}+a b_1^{21}+a^8b_1^{20}+b_1^{20}\cr
&+a^7 b_1^{19}+a^5 b_1^{19}+a^3b_1^{19}+a b_1^{19}+a^8 b_1^{18}+a^4 b_1^{18}+a^2
b_1^{18}+a^9 b_1^{17}+a^5 b_1^{17}+a^6 b_1^{16}\cr
&+a^4b_1^{16}) k^7+(a b_1^{21}+a^4 b_1^{20}+a^2b_1^{20}+b_1^{20}+a^7 b_1^{19}+a^{10} b_1^{18}+a^6b_1^{18}+a^4 b_1^{18}+b_1^{18}\cr
&+a^7 b_1^{17}+a^3b_1^{17}+a^{10} b_1^{16}+a^6 b_1^{16}+a^{11}
b_1^{15}+a^9 b_1^{15}+a^7 b_1^{15}+a^3 b_1^{15}+a^8b_1^{14}+a^6 b_1^{14}) k^6\cr
&+(a b_1^{19}+a^2b_1^{18}+b_1^{18}+a^7 b_1^{17}+a^5 b_1^{17}+a^{12}
b_1^{16}+a^8 b_1^{16}+a^6 b_1^{16}+a^4 b_1^{16}+a^2b_1^{16}+b_1^{16}\cr
&+a^{11} b_1^{15}+a^7 b_1^{15}+a^3b_1^{15}+a b_1^{15}+a^{12} b_1^{14}+a^{10}
b_1^{14}+a^6 b_1^{14}+a^4 b_1^{14}+a^2b_1^{14}+a^{13} b_1^{13}\cr
&+a^7 b_1^{13}+a^{10}b_1^{12}+a^8 b_1^{12}+a^4 b_1^{12}+a^5 b_1^{11}+a
b_1^{11}+a^6 b_1^{10}+a^2 b_1^{10}) k^5+(a^5b_1^{17}+a^3 b_1^{17}\cr
&+a^6 b_1^{16}+a^4 b_1^{16}+a^2b_1^{16}+a^{11} b_1^{15}+a^9 b_1^{15}+a^7 b_1^{15}+ab_1^{15}+a^{14} b_1^{14}+a^{10} b_1^{14}+a^8b_1^{14}\cr
&+a^6 b_1^{14}+a^2 b_1^{14}+a^9 b_1^{13}+a^3b_1^{13}+a b_1^{13}+a^{14} b_1^{12}+a^{12}b_1^{12}+a^8 b_1^{12}+a^2 b_1^{12}+a^{15}b_1^{11}\cr
&+a^{13} b_1^{11}+a^7 b_1^{11}+a^5b_1^{11}+a^{12} b_1^{10}+a^{10} b_1^{10}+a^8
b_1^{10}+a^6 b_1^{10}+a^2 b_1^{10}+a^8 b_1^8+a^4b_1^8) k^4\cr
&+(b_1^9 a^{17}+b_1^{12}a^{16}+b_1^{10} a^{16}+b_1^{11} a^{15}+b_1^8
a^{14}+b_1^{11} a^{13}+b_1^9 a^{13}+b_1^{10}a^{12}+b_1^8 a^{12}+b_1^{12} a^{10}\cr
&+b_1^{13}a^9+b_1^{11} a^9+b_1^{10} a^8+b_1^{13}a^7+b_1^9 a^7+b_1^{10} a^6+b_1^8 a^6+b_1^{15}a^5+b_1^{16} a^4+b_1^{10} a^4+b_1^8 a^4\cr
&+b_1^{15}a^3+b_1^{13} a^3+b_1^9 a^3+b_1^{14} a^2+b_1^{12}
a^2+b_1^{10} a^2+b_1^{13} a+b_1^{11} a)k^3+(b_1^7 a^{19}+b_1^{10} a^{18}\cr
&+b_1^8a^{18}+b_1^7 a^{17}+b_1^6 a^{16}+b_1^{11}a^{15}+b_1^9 a^{15}+b_1^7 a^{15}+b_1^{10}a^{14}+b_1^8 a^{14}+b_1^6 a^{14}+b_1^9 a^{13}\cr
&+b_1^8a^{12}+b_1^{11} a^{11}+b_1^7 a^{11}+b_1^{12}a^{10}+b_1^{10} a^{10}+b_1^8 a^{10}+b_1^{13}a^9+b_1^{11} a^9+b_1^9 a^9+b_1^7 a^9\cr
&+b_1^{12}a^8+b_1^6 a^8+b_1^{13} a^7+b_1^{11} a^7+b_1^9a^7+b_1^7 a^7+b_1^{14} a^6+b_1^{12} a^6+b_1^8a^6+b_1^6 a^6+b_1^{12} a^4\cr
&+b_1^8 a^4+b_1^{13}a^3+b_1^{11} a^3+b_1^{14} a^2+b_1^{10} a^2+b_1^{13}
a+b_1^{11} a) k^2+(b_1^2 a^{22}+b_1^3a^{21}+b_1^4 a^{20}\cr
&+b_1^5 a^{19}+b_1^6 a^{18}+b_1^4a^{18}+b_1^2 a^{18}+b_1^3 a^{17}+b_1^8 a^{16}+b_1^4a^{16}+b_1^9 a^{15}+b_1^5 a^{15}+b_1^{10}a^{14}\cr
&+b_1^2 a^{14}+b_1^9 a^{13}+b_1^3 a^{13}+b_1^8a^{12}+b_1^4 a^{12}+b_1^7 a^{11}+b_1^5 a^{11}+b_1^8a^{10}+b_1^6 a^{10}+b_1^4 a^{10}\cr
&+b_1^2 a^{10}+b_1^7a^9+b_1^3 a^9+b_1^{12} a^8+b_1^{10} a^8+b_1^4
a^8+b_1^{11} a^7+b_1^7 a^7+b_1^5 a^7+b_1^{10}a^6+b_1^8 a^6\cr
&+b_1^9 a^5+b_1^7 a^5+b_1^{12}
a^4+b_1^{10} a^4+b_1^{11} a^3+b_1^9 a^3) k.
\end{align*}
In \eqref{2.26.1} -- \eqref{2.28}, 
\begin{align*}
\tag{A6}\label{A6}
S_1=\,& b_1^8 k^3+b_1^4 k+b_1^3+b_1^2 k+b_1^2+b_1+1,\\
\tag{A7}\label{A7}
S_2=\,& b_1^{16} k^{15}+b_1^{16} k^{14}+b_1^{16} k^{13}+b_1^{16}
k^{12}+b_1^{14} k^{12}+b_1^{14} k^{10}+b_1^{12}k^{12}\\
&+b_1^{10} k^{12}+b_1^{10} k^{11}+b_1^{10}k^{10}+b_1^{10} k^8+b_1^{10} k^6+b_1^8 k^{12}+b_1^8k^{11}\cr
&+b_1^8 k^{10}+b_1^8 k^9+b_1^8 k^5+b_1^6
k^7+b_1^6 k^5+b_1^6 k^4+b_1^4 k^4+b_1^4k^3\cr
&+b_1^2 k^4+b_1^2 k^3+b_1^2 k^2+k^6+k^5+k+1,\cr
\tag{A8}\label{A8}
S_3=\,& b_1^{28} k^9+b_1^{26} k^7+b_1^{24} k^9+b_1^{24}k^8+b_1^{24} k^7+b_1^{24} k^6+b_1^{22} k^6+b_1^{20}k^{11}\\
&+b_1^{20} k^{10}+b_1^{20} k^8+b_1^{20}k^7+b_1^{20} k^5+b_1^{18} k^7+b_1^{18} k^6+b_1^{18}k^5\cr
&+b_1^{18} k^4+b_1^{16} k^9+b_1^{16} k^8+b_1^{16}k^7+b_1^{16} k^5+b_1^{16} k^3+b_1^{14} k^8\cr
&+b_1^{14}k^7+b_1^{14} k^6+b_1^{14} k^4+b_1^{12} k^8+b_1^{12}
k^3+b_1^{12} k^2+b_1^{10} k^6\cr
&+b_1^8 k^6+b_1^8k^4+b_1^8+b_1^4 k^2+1, \cr
\tag{A9}\label{A9}
T_1=\,&  a^8+a^6+a^5+a^2 k+a^2+a+1,\\
\tag{A10}\label{A10}
T_2=\,& a^{16} k^3+a^{16} k^2+a^{16} k+a^{16}+a^{14} k^2+a^{14}+a^{12}+a^{10}
k^3+a^8 k^4\\
&+a^6+a^4 k^2+a^4+a^2 k^3+a^2 k^2+a^2 k+k^4+k^3+k^2+k+1,\cr
\tag{A11}\label{A11}
T_3=\,& a^{26} k^2+a^{22} k^4+a^{22} k^2+a^{20} k^4+a^{20} k^3+a^{20} k^2+a^{18}
k^8+a^{18} k^7\\
&+a^{18} k^6+a^{18} k^4+a^{16} k^6+a^{16} k^3+a^{16}
k^2+a^{14} k^7+a^{14} k^5\cr
&+a^{14} k^4+a^{14} k^2+a^{12} k^8+a^{12}
k^5+a^{12} k^3+a^{12} k^2+a^{12} k\cr
&+a^{10} k^9+a^{10} k^8+a^{10} k^7+a^{10}
k^5+a^{10} k^2+a^{10} k+a^{10}+a^8 k^8\cr
&+a^8 k^7+a^8 k^4+a^6 k^5+a^6
k^3+a^6 k^2+a^4 k^8+a^4 k^7+a^4 k^6\cr
&+a^4 k^4+a^2 k^9+a^2 k^5+a^2 k^4+k^8.
\end{align*}
In \eqref{h1T1} -- \eqref{h3T1},
\begin{align*}
\tag{A12}\label{A12}
&h_1^*=a^{33}+a^{31}+a^{29}+a^{27}+(a^{61}+a^{60}+a^{59}+a^{58}+a^{57}+a^{56}+a^{53}+a^{52}+a^{51}\\
&+a^{50}+a^{49}+a^{48}+a^{37}+a^{36}+a^{35}+a^{34}+a^{33}+a^{32}
+a^{29}+a^{28}+a^{27}+a^{26}+a^{25}\cr
&+a^{24}+a^{21}+a^{20}+a^{19}+a^{18}+a^{17}+a^{16}+a^5+a^4+a^3+a^2+a+1)
b_1^{19}+(a^{66}\cr
&+a^{65}+a^{64}+a^{62}+a^{60}+a^{59}+a^{58}+a^{56}+a^{55}+a^{49}+a^{48}+a^{47}+a^{45}+a^{36}+a^{33}\cr
&+a^{32}+a^{28}+a^{27}+a^{24}+a^{23}+a^{22}+a^{18}+a^{17}+a^{16}+a^{15}+a^{14}+a^{13}+a^{10}+a^6\cr
&+a^4)b_1^{18}+(a^{67}+a^{66}+a^{60}+a^{59}+a^{58}+a^{57}+a^{55}+a^{54}+a
^{52}+a^{51}+a^{50}+a^{49}\cr
&+a^{46}+a^{45}+a^{44}+a^{40}+a^{38}+a^{33}+a^{32}+a^
{29}+a^{26}+a^{25}+a^{24}+a^{20}+a^{19}+a^{14}\cr
&+a^{11}+a^7+a^3+1)b_1^{17}+(a^{65}+a^{63}+a^{61}+a^{60}+a^{58}+a^{57}+a^{55}+a^{54}+a^{53}+a^{52}\cr
&+a^{50}+a^{48}+a^{46}+a^{41}+a^{38}+a^{36}+a^{32}+a^{30}+a^{28}+a^
{27}+a^{26}+a^{25}+a^{23}+a^{22}\cr
&+a^{15}+a^{13}+a^{11}+a^9+a^6+a^5+a^2+a)b_1^{16}+(a^{58}+a^{55}+a^{54}+a^{52}+a^{51}+a^{50}\cr
&+a^{49}+a^{47}+a^{46}+a^{44}+a^{42}+a^{41}+a^{38}+a^{36}+a^{33}+a^{32}+a^{30}+a^{28}+a^{27}+a^{26}\cr
&+a^{20}+a^{19}+a^{16}+a^{15}+a^{11}+a^{10}+a^9+a^7+a^4+a^2)
b_1^{15}+(a^{56}+a^{55}+a^{53}+a^{52}\cr
&+a^{48}+a^{46}+a^{45}+a^{44}+a^{42}+a^{35}+a^{33}+a^{31}+a^{30}+a^{27}+a^{26}+a^{25}+a^{24}+a^{23}\cr
&+a^{22}+a^{21}+a^{20}+a^{18}+a^{17}+a^{15}+a^{14}+a^{13}+a^{11}+a^{10}+a^9+a^8+a^4+a^3)b_1^{14}\cr
&+(a^{54}+a^{52}+a^{51}+a^{50}+a^{49}+a^{48}+a^{46}+a^{45}+a
^{40}+a^{38}+a^{36}+a^{35}+a^{27}+a^{26}\cr
&+a^{22}+a^{21}+a^{19}+a^{18}+a^{17}+a^{16}+a^{12}+a^{10}+a^8+a^4)
b_1^{13}+(a^{52}+a^{51}+a^{50}\cr
&+a^{48}+a^{47}+a^{41}+a^{37}+a^{36}+a^{35}+a^{34}+a^{33}+a^{30}+a^{29}+a^{28}+a^{27}+a^{24}+a^{13}\cr
&+a^{12}+a^{11}+a^7+a^6+a^5)b_1^{12}+(a^{50}+a^{47}+a^{44}+a^{42}+a^{40}+a^{38}+a^{37}+a^{35}\cr
&+a^{34}+a^{33}+a^{27}+a^{25}+a^{23}+a^{22}+a^{21}+a^{17}+a^{16}+a^{14}+a^{12}+a^{10}+a^9+a^6)b_1^{11}\cr
&+(a^{48}+a^{47}+a^{45}+a^{43}+a^{42}+a^{33}+a^{32}+a^{31}+a
^{29}+a^{27}+a^{26}+a^{20}+a^{17}+a^{16}\cr
&+a^{15}+a^{12}+a^8+a^7)b_1^{10}+(a^{46}+a^{44}+a^{42}+a^{41}+a^{39}+a^{38}+a^{36}+a^{35}+a^{34}\cr
&+a^{32}+a^{31}+a^{30}+a^{28}+a^{27}+a^{26}+a^{24}+a^{23}+a^{21}+a^{20}+a^
{19}+a^{18}+a^{16}+a^{15}\cr
&+a^{14}+a^{12}+a^{11})b_1^9+(a^{44}+a^{43}+a^{42}+a^{41}+a^{39}+a^{38}+a^{37}+a^{36}+a^{34}+a^{32}\cr
&+a^{31}+a^{29}+a^{28}+a^{27}+a^{26}+a^{24}+a^{23}+a^{22}+a^{20}+a^{
17}+a^{16}+a^{15}+a^{13}+a^{11})b_1^8\cr
&+(a^{42}+a^{39}+a^{38}+a^{36}+a^{35}+a^{34}+a^{33}+a^{31}+a^{
30}+a^{26}+a^{25}+a^{24}+a^{17})b_1^7\cr
&+(a^{40}+a^{39}+a^{37}+a^{36}+a^{32}+a^{30}+a^{29}+a^{28}+a^{
26}+a^{20}+a^{19}+a^{18}+a^{17})b_1^6\cr
&+(a^{38}+a^{36}+a^{35}+a^{34}+a^{33}+a^{32}+a^{30}+a^{29}+a^{25}+a^{23}
+a^{22}+a^{21}+a^{20})b_1^5\cr
&+(a^{36}+a^{35}+a^{34}+a^{32}+a^{31}+a^{25}+a^{23}+a^{22}) 
b_1^4+(a^{34}+a^{31}+a^{29}+a^{27}+a^{25}\cr
&+a^{22})b_1^3+(a^{34}+a^{33}+a^{31}+a^{30}+a^{29}+a^{28}+a^{24}+a^{23}
) b_1^2+(a^{35}+a^{34}+a^{30}+a^{27}) b_1, 
\end{align*}
\begin{align*}
\tag{A13}\label{A13}
&h_2^*= a^{35}+a^{33}+a^{31}+a^{29}+(a^{73}+a^{72}+a^{69}+a^{68}+a^{67}+a^{66}+a^{65}+a^{64}+a^{57}\\
&+a^{56}+a^{53}+a^{52}+a^{51}+a^{50}+a^{45}+a^{44}+a^{43}+a^{42}
+a^{41}+a^{40}+a^{25}+a^{24}+a^{21}\cr
&+a^{20}+a^{19}+a^{18}+a^{13}+a^{12}+a^{11}+
a^{10}+a^5+a^4+a^3+a^2+a+1)b_1^{19}+(a^{71}\cr
&+a^{69}+a^{65}+a^{55}+a^{53}+a^{49}+a^{47}+a^{45}+a
^{41}+a^{23}+a^{21}+a^{17}+a^{15}+a^{13}+a^9\cr
&+a^7+a^5+a)b_1^{18}+(a^{69}+a^{68}+a^{64}+a^{63}+a^{60}+a^{59}+a^{56}+a^{55}+a
^{53}+a^{51}+a^{45}\cr
&+a^{44}+a^{40}+a^{39}+a^{36}+a^{35}+a^{32}+a^{31}+a^{28}+a^
{27}+a^{24}+a^{23}+a^{21}+a^{19}+a^{13}\cr
&+a^{12}+a^8+a^7+a^5+a^3)b_1^{17}+(a^{67}+a^{65}+a^{63}+a^{60}+a^{59}+a^{57}+a^{56}+a^{54}+a^{52}\cr
&+a^{51}+a^{50}+a^{49}+a^{47}+a^{44}+a^{42}+a^{41}+a^{40}+a^{39}+a^{37}+a^
{32}+a^{29}+a^{23}+a^{22}\cr
&+a^{21}+a^{18}+a^{16}+a^{13}+a^{11}+a^{10}+a^4)b_1^{16}+(a^{59}+a^{57}+a^{55}+a^{54}+a^{53}+a^{51}\cr
&+a^{50}+a^{48}+a^{47}+a^{45}+a^{40}+a^{37}+a^{36}+a^{35}+a^{34}+a^{30}+a^{29}+a^{27}+a^{24}+a^{22}\cr
&+a^{21}+a^{20}+a^{16}+a^{15}+a^{14}+a^9+a^7+a^5)b_1^{15}+(a^{57}+a^{56}+a^{55}+a^{54}+a^{49}+a^{47}\cr
&+a^{43}+a^{41}+a^{39}+a^{37}+a^{33}+a^{32}+a^{30}+a^{25}+a^{23}+a^{21}+a^{16}+a^{15}+a^{14}+a^{11}\cr
&+a^9+a^8+a^7+a^6)b_1^{14}+(a^{55}+a^{47}+a^{46}+a^{45}+a^{44}+a^{41}+a^{39}+a^{36}+a^{35}+a^{33}\cr
&+a^{32}+a^{31}+a^{29}+a^{26}+a^{25}+a^{23}+a^{20}+a^{19}+a^{17}+a^{16}+a^{14}+a^{12}+a^{10}+a^7)b_1^{13}\cr
&+(a^{53}+a^{52}+a^{49}+a^{47}+a^{46}+a^{42}+a^{41}+a^{40}+a
^{38}+a^{36}+a^{34}+a^{30}+a^{28}+a^{27}\cr
&+a^{24}+a^{21}+a^{20}+a^{17}+a^{16}+a^{15}+a^{14}+a^{11}+a^9+a^8)b_1^{12}+(a^{51}+a^{48}+a^{46}\cr
&+a^{45}+a^{44}+a^{41}+a^{37}+a^{35}+a^{34}+a^{32}+a^{30}+a^{27}+a^{25}+a^{19}+a^{16}+a^{12}+a^{10}\cr
&+a^8)b_1^{11}+(a^{49}+a^{48}+a^{46}+a^{44}+a^{43}+a^{42}+a^{41}+a^{39}+a
^{37}+a^{33}+a^{27}+a^{26}\cr
&+a^{22}+a^{21}+a^{19}+a^{18}+a^{16}+a^{15}+a^{12}+a^{11}+a^{10}+a^9)
b_1^{10}+(a^{47}+a^{45}+a^{44}\cr
&+a^{43}+a^{40}+a^{39}+a^{37}+a^{35}+a^{31}+a^{30}+a^{29}+a^{27}+a^{20}+a^{19}+a^{16}+a^{15}+a^{14}\cr
&+a^{13})b_1^9+(a^{45}+a^{44}+a^{43}+a^{41}+a^{40}+a^{38}+a^{36}+a^{35}+a^{
34}+a^{33}+a^{31}+a^{30}\cr
&+a^{25}+a^{24}+a^{18}+a^{17}+a^{15}+a^{13})b_1^8+(a^{43}+a^{41}+a^{39}+a^{38}+a^{37}+a^{35}+a^{34}\cr
&+a^{32}+a^{31}+a^{29}+a^{28}+a^{22}+a^{21}+a^{18})
b_1^7+(a^{41}+a^{40}+a^{39}+a^{38}+a^{33}+a^{31}\cr
&+a^{27}+a^{25})b_1^6+(a^{39}+a^{31}+a^{30}+a^{29}+a^{28}+a^{25}+a^{22})
b_1^5+(a^{37}+a^{36}+a^{33}+a^{31}\cr
&+a^{30}+a^{26}+a^{25}+a^{24})b_1^4+(a^{41}+a^{40}+a^{37}+a^{36}+a^{34}+a^{33}+a^{30}+a^{29}+a^{28}\cr
&+a^{24})b_1^3+(a^{39}+a^{37}+a^{32}+a^{30}+a^{28}+a^{27}+a^{26}+a^{25}
) b_1^2+(a^{37}+a^{36}+a^{32}+a^{29}) b_1,  
\end{align*}
\begin{align*}
\tag{A14}\label{A14}
&h_3^*= a^{36}+a^{28}+(a^{68}+a^{67}+a^{66}+a^{64}+a^{62}+a^{60}+a^{58}+a^{56}+a^{54
}+a^{51}+a^{44}\\
&+a^{43}+a^{42}+a^{40}+a^{38}+a^{36}+a^{34}+a^{32}+a^{30}+a^{28}
+a^{26}+a^{24}+a^{22}+a^{19}+a^{12}\cr
&+a^{11}+a^{10}+a^8+a^6+a^3)b_1^{23}+(a^{67}+a^{65}+a^{61}+a^{59}+a^{58}+a^{57}+a^{53}+a^{50}+a^{43}\cr
&+a^{41}+a^{37}+a^{35}+a^{34}+a^{33}+a^{29}+a^{27}+a^{26}+a^{25}+a^{21}+a^
{18}+a^{11}+a^9+a^5\cr
&+a^2)b_1^{22}+(a^{68}+a^{67}+a^{66}+a^{65}+a^{64}+a^{61}+a^{60}+a^{59}+a
^{57}+a^{56}+a^{55}+a^{54}\cr
&+a^{53}+a^{52}+a^{51}+a^{50}+a^{48}+a^{47}+a^{45}+a^{44}+a^{43}+a^{38}+a^{37}+a^{36}+a^{34}+a^{32}\cr
&+a^{30}+a^{29}+a^{28}+a^{27}+a^{26}+a^{25}+a^{24}+a^{23}+a^{21}+a^{20}+a^{19}+a^{16}+a^{15}+a^{14}\cr
&+a^{13}+a^{12}+a^{11}+a^{10}+a^5+a^3+a^2+a)
b_1^{21}+(a^{69}+a^{68}+a^{63}+a^{60}+a^{57}+a^{56}\cr
&+a^{55}+a^{50}+a
^{49}+a^{47}+a^{43}+a^{41}+a^{40}+a^{38}+a^{37}+a^{36}+a^{35}+a^{32}+a^{30}+a^
{28}\cr
&+a^{27}+a^{24}+a^{23}+a^{21}+a^{20}+a^{17}+a^{16}+a^{14}+a^8+a^6+a^5+a^4+a
^3+a^2)b_1^{20}\cr
&+(a^{68}+a^{66}+a^{65}+a^{64}+a^{63}+a^{62}+a^{59}+a^{58}+a
^{57}+a^{55}+a^{54}+a^{50}+a^{48}+a^{47}\cr
&+a^{45}+a^{42}+a^{40}+a^{37}+a^{36}+a^{33}+a^{31}+a^{29}+a^{28}+a^{25}+a^{24}+a^{21}+a^{20}+a^{15}\cr
&+a^{14}+a^{12}+a^{11}+a^7+a^6+a^4)b_1^{19}+(a^{69}+a^{68}+a^{67}+a^{66}+a^{65}+a^{63}+a^{54}+a^{53}\cr
&+a^{52}+a^{51}+a^{50}+a^{49}+a^{48}+a^{46}+a^{42}+a^{41}+a^{39}+a^{38}+a^{30}+a^{29}+a^{26}+a^{21}\cr
&+a^{20}+a^{19}+a^{15}+a^{13}+a^9+a^7+a^5+a^4+a^3+1)
b_1^{18}+(a^{70}+a^{69}+a^{68}+a^{67}\cr
&+a^{66}+a^{64}+a^{62}+a^{61}+a^{60}+a^{59}+a^{58}+a^{55}+a^{54}+a^{52}+a^{48}+a^{47}+a^{45}+a^{44}\cr
&+a^{42}+a^{41}+a^{40}+a^{38}+a^{37}+a^{36}+a^{35}+a^{33}+a^{32}+a^{30}+a^{28}+a^{27}+a^{26}+a^{24}\cr
&+a^{23}+a^{22}+a^{20}+a^{19}+a^{18}+a^{17}+a^{15}+a^{14}+a^9+a^8+a^6
+a^4+a^3+a^2)b_1^{17}\cr
&+(a^{68}+a^{64}+a^{63}+a^{62}+a^{60}+a^{58}+a^{53}+a^{52}+a
^{50}+a^{46}+a^{42}+a^{40}+a^{39}+a^{36}\cr
&+a^{35}+a^{29}+a^{26}+a^{25}+a^{24}+a^{23}+a^{21}+a^{20}+a^{19}+a^{18}+a^{17}+a^{15}+a^{14}+a^{13}\cr
&+a^{12}+a^{11}+a^{10}+a^8+a^6+a^3)b_1^{16}+(a^{62}+a^{59}+a^{58}+a^{54}+a^{52}+a^{48}+a^{47}+a^{42}\cr
&+a^{40}+a^{39}+a^{38}+a^{33}+a^{31}+a^{30}+a^{29}+a^{26}+a^{20}+a^{18}+a^{17}+a^{14}+a^{13}+a^{12}\cr
&+a^{11}+a^7+a^6+a^4)b_1^{15}+(a^{60}+a^{59}+a^{55}+a^{54}+a^{53}+a^{48}+a^{44}+a^{42}+a^{41}+a^{40}\cr
&+a^{39}+a^{38}+a^{37}+a^{34}+a^{32}+a^{30}+a^{29}+a^{28}+a^{27}+a^
{26}+a^{23}+a^{21}+a^{18}+a^{16}\cr
&+a^{13}+a^{12}+a^9+a^7+a^6+a^5)b_1^{14}+(a^{58}+a^{56}+a^{52}+a^{50}+a^{49}+a^{47}+a^{46}+a^{44}\cr
&+a^{43}+a^{39}+a^{38}+a^{33}+a^{32}+a^{28}+a^{26}+a^{24}+a^{22}+a^{21}+a^{19}+a^{18}+a^{16}+a^{13}\cr
&+a^{12}+a^6)b_1^{13}+(a^{56}+a^{55}+a^{54}+a^{49}+a^{47}+a^{44}+a^{42}+a^{41}+a
^{39}+a^{38}+a^{37}\cr
&+a^{36}+a^{34}+a^{33}+a^{30}+a^{29}+a^{28}+a^{27}+a^{26}+a^
{22}+a^{18}+a^{17}+a^{15}+a^{14}+a^{12}\cr
&+a^{10}+a^8+a^7)b_1^{12}+(a^{54}+a^{51}+a^{49}+a^{48}+a^{46}+a^{44}+a^{42}+a^{40}+a
^{30}+a^{26}\cr
&+a^{25}+a^{24}+a^{21}+a^{19}+a^{16}+a^{14}+a^{12}+a^{11}+a^{10}+a^8)
b_1^{11}+(a^{52}+a^{51}+a^{48}\cr
&+a^{47}+a^{46}+a^{45}+a^{43}+a^{40}+a^{37}+a^{35}+a^{33}+a^{32}+a^{31}+a^{25}+a^{22}+a^{21}+a^{20}\cr
&+a^{19}+a^{18}+a^{16}+a^{14}+a^{12}+a^{10}+a^9+a^8)
b_1^{10}+(a^{50}+a^{48}+a^{42}+a^{41}+a^{40}\cr
&+a^{38}+a^{33}+a^{32}+a^{31}+a^{30}+a^{29}+a^{28}+a^{27}+a^{24}+a^{23}+a^{22}+a^{21}+a^{19}+a^{16}\cr
&+a^{13}+a^{12}+a^{11})b_1^9+(a^{48}+a^{47}+a^{46}+a^{42}+a^{39}+a^{36}+a^{33}+a^{32}+a^{29}+a^{27}\cr
&+a^{26}+a^{23}+a^{20}+a^{19}+a^{15}+a^{14}+a^{12})
b_1^8+(a^{46}+a^{43}+a^{42}+a^{38}+a^{37}+a^{36}\cr
&+a^{34}+a^{27}+a^{25}+a^{23}+a^{21}+a^{19}+a^{18}+a^{17})
b_1^7+(a^{44}+a^{43}+a^{39}+a^{38}+a^{37}\cr
&+a^{32}+a^{30}+a^{29}+a^{27}+a^{23}+a^{18})
b_1^6+(a^{42}+a^{40}+a^{35}+a^{33}+a^{30}+a^{29}+a^{28}\cr
&+a^{23})b_1^5+(a^{40}+a^{39}+a^{38}+a^{37}+a^{36}+a^{33}+a^{31}+a^{29}+a^{
25}+a^{22})b_1^4+(a^{38}\cr
&+a^{36}+a^{35}+a^{34}+a^{31}+a^{29}+a^{25}+a^{24})
b_1^3+(a^{37}+a^{34}+a^{33}+a^{32}+a^{30}+a^{29}\cr
&+a^{25}+a^{24})b_1^2+(a^{38}+a^{37}+a^{36}+a^{35}+a^{30}+a^{29}+a^{28}+a^{27}
) b_1. 
\end{align*}

\end{document}